\documentclass[10pt]{amsart}
\usepackage{geometry}                % See geometry.pdf to learn the layout options. There are lots.
\usepackage{graphicx}
\usepackage{amssymb}
\usepackage{epstopdf}
\usepackage{tikz}
\newtheorem{Thm}{Theorem}

\newtheorem{Def}{Definition}
\newtheorem{Cor}{Corollary}
\address{Institute for Mathematical Behavioral Sciences, University of California, Irvine; Irvine CA 92697-5100}
\email{dsaari@uci.edu} 
\title[Symmetry structures]{Inherent Symmetries of graphs,  paths, and Traveling Salesperson Problems}
\author[D. Saari ]{Donald G. Saari}   
\thanks{My thanks to  George Hazelrigg for our several  discussions.  This work is part of a National Science Foundation project under NSF Award Number CMMI-1923164.}
\begin{document}
\maketitle
%\section{}
%\subsection{}
%\begin{abstract}   
 \begin{abstract} Without imposing   restrictions   on a weighted  graph's  arc lengths,  symmetry structures cannot be expected.  But, they exist.  To find them, the  graphs are decomposed  into a component that dictates all  closed path properties (e.g.,  shortest and longest   paths),  and a superfluous   component  that can be removed.  The simpler remaining graph exposes inherent symmetry structures that form  the basis for all  closed  path properties.  For certain asymmetric problems, the symmetry is that   of three-cycles; for the general undirected setting it is a type of  four-cycles; for general directed problems with asymmetric costs, it is a product of three and four cycles.  Everything extends  immediately to incomplete graphs. 

 \end{abstract}
 \medskip
 %\begin{keywords}   Traveling salesperson problem, optimal paths, greedy algorithm,  closed path symmetries.\end{keywords}
%\keywords{Traveling salesperson problem, optimal paths, greedy algorithm,  closed path symmetries}

%\newpage

\section{Introduction}  \label{sect: intro}   With  
  applications   ranging from the design of microchips to the positioning of telescopes \cite{cook}, understanding properties of a weighted  graph's paths and  closed paths  has attained  importance beyond  mathematics.    
What complicates this analysis is that   transmission costs  between nodes typically   include  factors  
 that differ  
 from what is needed to    determine path properties; e.g., they may reflect  
the problem's topography or congestions of various types.    
These portions  of a graph's entries  add nothing to the analysis,  but  they contribute to the    complexity of   these concerns.

 The mathematical structure of graphs       developed here 
  separates a complete  weighted  graph into two  unique  components.  The first,  
 with  best possible  reduced    degrees of freedom, has     all of the information   needed to develop the particulars of paths and  closed paths.     The remaining component is dismissed because it adds no value; it  
just  complicates both the    analysis of paths and the behavior of algorithms.      This approach  
   resembles  (and  is  motivated by)   a   game theory decomposition \cite{jessie} where one game  component   has       only (and all)  information needed to find all pure and mixed Nash strategic properties; another component  captures coordination,  cooperation, etc.

Three classes of graphs are examined.  The first  is a directed, asymmetric  setting measuring 
  differences from the average cost between vertices.  The second and third are, respectively, 
 the standard undirected  symmetric cost     and directed asymmetric cost settings.  Symmetry structures for these  classes differ; e.g., the  symmetry structures for
the   excess cost  graphs (Sect.~2)  are three-cycles.  Symmetries for the standard symmetric case (Sect.~3)  are a form of four-cycles.  Symmetries for the general  asymmetric costs (Sect.~4)  are  a product of three and four cycles.

To describe the basic theme in terms of the first class, it turns out that 
 these graphs can be embedded in
  the  
space of   asymmetric   paired comparisons.  A  ``decision theory"  decomposition  divides this space   
 into a linear subspace characterized by 
  a  strong form of transitivity and its normal bundle  consisting  of cycles \cite{tandd}.    
Voting  methods seek     linear orders,    
 so  the  cyclic components    create    
   complexities and paradoxical outcomes   (e.g., Arrow's Impossibility Theorem \cite{arrow}).
Projecting the data   to the transitive subspace eliminates these difficulties and  simplifies the analysis \cite{tandd}.   But 
 cycles, not linear orders, are central for closed paths,   
  so   in this setting 
    the   transitive components    are what 
        obscure  
   the analysis.    Projecting the   data   to the cyclic subspace lowers the  
 degrees of freedom,  
  removes  trouble-causing  components, and uncovers  the system's inherent three-cycle symmetry. 
Here,    a $A, B, C$ cycle's costs  of going from $A$ to $B$,    $B$ to $C$, and   $C$ to $A$ are identical. 
  
In general,  each   class  of   graphs is decomposed  
 into a component characterized by a  behavior that masks the  closed path structures   and a component that has only (and all) of the  relevant closed path information.     Everything extends to incomplete graphs. Most proofs are in Sect.~\ref{sect: proofs}.

\section{Asymmetric excess costs}

Reimbursing a salesperson for the average cost of traveling between   cities creates an incentive to find    routes  with 
 below  average costs.     For notation,  
 it takes 40 minutes to walk  from   home, H, to campus, C; returning uphill requires 50 minutes, so the average is 45.  
  The   ``excess cost function"    registers    differences from the average  where $ C \stackrel{5}{\longrightarrow} H \textrm{ represents both  } C \stackrel{5}{\longrightarrow} H \textrm{ and }H \stackrel{-5}{\longrightarrow} C. $

  Graphs in the space of  asymmetric weighted, n-vertex graphs (with no loops) considered here, $\mathbb G^n_A$,  are complete (i.e., each pair of vertices is connected with  paths) and \begin{equation}\label{def: dij=-dji} V_j  \stackrel{x}{\longrightarrow} V_k \textrm{ represents both  }    V_j  \stackrel{x}{\longrightarrow} V_k \textrm{ and }V_k \stackrel{-x}{\longrightarrow} V_j.\end{equation} To simplify the  
   graphs,  
   only an arc's positive cost direction  need be represented; this is because moving counter to an arrow represents a ``below average" cost (Eq.~\ref{def: dij=-dji}). With this choice, $V_2$ in Fig.~1a  is a ``source" as all positive value directions point away; it is  a ``sink"  with the  negative value directions. Conversely, $V_4$ is a sink with positive value directions and a source for negative value directions.  Subscripts $A$ and $S$ indicate, respectively, the asymmetric and symmetric cases.  
   
   \begin{tikzpicture}[xscale=0.5, yscale=0.52]

\node[below] at (1, 0) {$V_1$};
 
\node[below] at (4, 0) {$V_2$};

\node[right] at (5, 2) {$V_3$};
 
\node[above] at (2.5, 3) {$V_4$};
 
\node[left] at (0, 2) {$V_5$};
\draw[<-] (1.1,  0) -- (3.8, 0); %ab
\draw [<-] (4.9, 1.9) -- (1.15, .15);  %ca
\draw [->]  (1, 0) -- (2.4, 2.9); %ad
\draw [->] (.9, .05) -- (.1, 1.9); %ae
\draw [<-]  (0.15, 2) -- (3.8, .15);  %be
 \draw [<-] (2.55, 3) -- (3.9, .15);  %db
 \draw [->]  (4, 0) -- (4.9, 1.8); %bc
\draw [->]  (4.9, 2.1) -- (2.75, 3); %cd
\draw [->] (4.85, 2) -- (0.45, 2); %ce
\draw [->] (0.2, 2.1) -- (2.3, 2.95); %ed

\node[below] at (2.4, 0) {$1$};  %ab
\node[right] at (4.5, 1) {$3$}; %bc
\node[left] at (0.5, 1) {$8$}; %ae
\node[above] at (3.75, 2.5) {$6$};  %cd
\node[above] at (1.25, 2.5) {$3$}; %ed

\node[right] at (2.54, .99) {\small 4}; %AC
 
\node[left] at (2.3, 1.75) {\small 14}; %ad
 
\node[right] at (1.56, .99) {\small 9};
 
\node[above] at (3.23, 1.3) {\small 7}; %bd
%\node[above] at (2.5, 2) {$0$}; %ec
\node[left] at (2.85, 2.1) {\small 1}; %ec

\node[below] at (2.5, -1) {{\bf a.} Graph $\mathcal G^5_{A}$}; 
\node[below] at (12.5, -2) {{\bf Figure 1.} Decomposition};
\node at (7.5, 2) {=};
\node at (17.5, 2) {+};

%%%%%%%%%%%%%%%

\node[below] at (11, 0) {$V_1$};
%\node at (24, 0) {$\bullet$};
\node[below] at (14, 0) {$V_2$};

%\node at (25, 2) {$\bullet$};
\node[right] at (15, 2) {$V_3$};
%\node at (22.5, 3) {$\bullet$};
\node[above] at (12.5, 3) {$V_4$};
%\node at (20, 2) {$\bullet$};
\node[left] at (10, 2) {$V_5$};
\draw[->] (11.1,  0) -- (13.8, 0); %ab
\draw [<-] (14.9, 1.9) -- (11.15, .15);  %ca
\draw [->]  (11, 0) -- (12.4, 2.9); %ad
\draw [->] (10.9, .05) -- (10.1, 1.9); %ae
\draw [<-]  (10.15, 2) -- (13.8, .15);  %be
 \draw [<-] (12.55, 3) -- (13.9, .15);  %db
 \draw [->]  (14, 0) -- (14.9, 1.8); %bc
\draw [->]  (14.9, 2.1) -- (12.75, 3); %cd
\draw [->] (14.85, 2) -- (10.45, 2); %ce
\draw [->] (10.2, 2.1) -- (12.3, 2.95); %ed

\node[below] at (12.4, 0) {$1$};  %ab
\node[right] at (14.5, 1) {$4$}; %bc
\node[left] at (10.5, 1) {$8$}; %ae
\node[above] at (13.75, 2.5) {$6$};  %cd
\node[above] at (11.25, 2.5) {$3$}; %ed

%\node[below] at (3.1, 1.3) {$0$}; %ac
\node[right] at (12.64, .99) {\small 5}; %AC
%\node[left] at (1.8, 1.75) {$2$}; %ad
\node[left] at (12.34, 1.75) {11}; %ad (2.3, 1.75
%\node[below] at (1.9, 1.3) {$1$}; %eb
\node[right] at (11.56, .99) {\small 7};
%\node[right] at (3.2, 1.75) {$0$}; %bd
\node[above] at (13.3, 1.3) {\small 10}; %bd
%\node[above] at (2.5, 2) {$0$}; %ec
\node[left] at (12.9, 2.1) {\small 3}; %ec

\node[below] at (12.2, -1) {{\bf b.} Closed path independent; $\mathcal G^5_{A, cpi}$}; 

%%%%%%%%%%%%%%%%%%%%%
%ab=1, ac=-4, ad=-3, ae=-9, bc=-5, bd=-4, be=-10; cd=1, ce=-5, de=-6
\node[below] at (21, 0) {$V_1$};
%\node at (24, 0) {$\bullet$};
\node[below] at (24, 0) {$V_2$};

%\node at (25, 2) {$\bullet$};
\node[right] at (25, 2) {$V_3$};
%\node at (22.5, 3) {$\bullet$};
\node[above] at (22.5, 3) {$V_4$};
%\node at (20, 2) {$\bullet$};
\node[left] at (20, 2) {$V_5$};
\draw[<-] (21.1,  0) -- (23.8, 0); %ab
\draw [->] (24.9, 1.9) -- (21.15, .15);  %ca
\draw [->]  (21, 0) -- (22.4, 2.9); %ad
\draw [dashed] (20.9, .05) -- (20.1, 1.9); %ae
\draw [<-]  (20.15, 2) -- (23.8, .15);  %be
 \draw [->] (22.55, 3) -- (23.9, .15);  %db
 \draw [<-]  (24, 0) -- (24.9, 1.8); %bc
\draw [dashed]  (24.9, 2.1) -- (22.75, 3); %cd
\draw [<-] (24.85, 2) -- (20.45, 2); %ce
\draw [dashed] (20.2, 2.1) -- (22.3, 2.95); %ed

\node[below] at (22.4, 0) {$2$};  %ab
\node[right] at (24.5, 1) {$1$}; %bc
\node[left] at (20.5, 1) {$0$}; %ae
\node[above] at (23.75, 2.5) {$0$};  %cd
\node[above] at (21.25, 2.5) {$0$}; %ed

%\node[below] at (3.1, 1.3) {$0$}; %ac
\node[right] at (22.64, .99) {\small 1}; %AC
%\node[left] at (1.8, 1.75) {$2$}; %ad
\node[left] at (22.3, 1.75) {\small 3}; %ad
%\node[below] at (1.9, 1.3) {$1$}; %eb
\node[right] at (21.56, .99) {\small 2};
%\node[right] at (3.2, 1.75) {$0$}; %bd
\node[above] at (23.3, 1.3) {\small 3}; %bd
%\node[above] at (2.5, 2) {$0$}; %ec
\node[left] at (22.9, 2.1) {\small 2}; %ec

\node[below] at (23.25, -1) {{\bf c.} Cyclic component; $\mathcal G^5_{A, cyclic}$};

\end{tikzpicture}

 Figure 1 depicts  the general approach whereby graph $\mathcal G^5_A$ is uniquely  decomposed into   a ``closed path independent" component $\mathcal G^5_{A, cpi}$ (defined in Def.~\ref{def: path independent}) and a cyclic component $\mathcal G^5_{A, cyclic}$ to have $ \mathcal G^5_A= \mathcal G^5_{A, cpi} + \mathcal G^5_{A, cyclic}.$ The goal is to achieve  this decomposition  for any $\mathcal G^n_A \in \mathbb G^n_A$ to obtain     \begin{equation}\label{eq: first decomp}  \mathcal G^n_A = \mathcal G^n_{A, cpi} + \mathcal G^n_{A, cyclic}.\end{equation}

 \begin{Def}\label{def: path independent} Graph $\mathcal G^n_{A, cpi} \in \mathbb G^n_A$  is ``closed path independent" (cpi) iff all closed paths have length zero.
  A graph is strongly transitive iff  path lengths of a   triplet $\{V_i, V_j, V_k\}$ satisfy  
  \begin{equation}\label{eq: st} 
 V_i \stackrel{x}{\longrightarrow} V_j  \stackrel{y}{\longrightarrow} V_k = V_i \stackrel{z=x+y}{\longrightarrow} V_k.\end{equation}  \end{Def}

Both  Eq.~\ref{eq: st} paths start at $V_i$ and end at $V_k$, so  equality designates 
 equal path lengths.
This equation
 modifies  the concept of %
  ``strong transitivity"  developed for 
   decision  theory \cite{tandd}.

  \begin{Thm} \label{thm: st} A graph is strongly transitive iff it  is  
  cpi.
  Strongly transitive graphs (equivalently,   cpi graphs)  with $n$ vertices define a $(n-1)$-dimensional linear subspace $\mathbb{ST}_A^n \subset \mathbb G^n_A$.\end{Thm}

The Fig.~1b graph is 
   strongly transitive and  cpi.  
To check for strong transitivity,  
 select any triplet, say $\{V_1, V_3, V_5\},$ and determine whether  this triangle's leg lengths,  $V_1  \stackrel{5}{\longrightarrow} V_3,   V_3 \stackrel{3}{\longrightarrow}  V_5,$ and $V_1 \stackrel{8}{\longrightarrow}  V_5$,   satisfy the {\em triangle equality} Eq.~\ref{eq: st}, which they do.  To equate   strong transitivity with cpi,   
 reversing  
 $V_1 \stackrel{8}{\longrightarrow}  V_5$ defines  
    the   closed path $V_1  \stackrel{5}{\longrightarrow} V_3  
  \stackrel{3}{\longrightarrow}  V_5  
  \stackrel{-8}{\longrightarrow}  V_1$ with zero length.   

While the proof of  
 Thm.~\ref{thm: st} is in Sect.~\ref{sect: proofs}, proving  
  that $\mathbb {ST}_A^n$ is a linear subspace is a common exercise.  For the dimensionality assertion,  strong transitivity  ensures that   
 the  $V_i \to V_j$ arc  length  equals the  $V_i \stackrel{x}{\longrightarrow} V_1  \stackrel{y}{\longrightarrow} V_j$ length,   where the  $V_i \to V_j$  path is  diverted to pass
  through $V_1$.     
As all  arc  lengths for   $\mathcal G^n_{A, cpi}\in \mathbb {ST}^n_A$  are determined by    the    $\{V_1\to V_k\}_{k=2}^n$ arc lengths,  
    $\mathbb{ST}_A^n$ has dimension $(n-1)$.  
     
 A standard induction argument applied to Eq.~\ref{eq: st}  proves the following result.
 
 \begin{Cor}\label{cor: gen path length} For $\mathcal G^n_{A, cpi} \in \mathbb {ST}_A^n$,  any path starting at $V_i$ and ending at $V_j$ has length equal to the  $V_i\to V_j$ arc that  connects the endpoints.\end{Cor}
 
 A Fig.~1b example of Cor.~\ref{cor: gen path length}  
  is where the $4$ length of   $V_2  \stackrel{7}{\longrightarrow} V_5  \stackrel{3}{\longrightarrow} V_4  \stackrel{-10}{\longrightarrow} V_2  \stackrel{10}{\longrightarrow} V_6  \stackrel{-6}{\longrightarrow} V_3$, where vertices can be revisited,    equals the arc length connecting the endpoints $V_2  \stackrel{4}{\longrightarrow} V_3$.

\subsection{Cyclic Normal Bundle.}   
 As   
Thm.~\ref{thm: gen decomp} will assert, the  $\mathcal G^n_{A, cpi} \in \mathbb{ST}_A^n$ component of  $\mathcal G^n_A$  (see Eq.~\ref{eq: first decomp}) blurs  the  
 closed path properties of $\mathcal G^n_A$.  Thus,      path properties must  be   
  based on the structure of   $\mathbb{ST}_A^n$'s normal bundle.   The $\mathbb G^n_A$  and $\mathbb{ST}^n_A$ dimensions are  $n\choose2$ and   
  $(n-1)$,  so  $\mathbb{ST}^n_A$'s   normal subspace,  $\mathbb C^n_A$,  has dimension ${{n-1}\choose 2}$. As described next,    $\mathbb C^n_A$ consists of cyclic actions.
 \begin{Thm}\label{thm: cn} \cite{tandd}  For $n\ge 3$, the linear subspace orthogonal to $\mathbb{ST}_A^n$, $\mathbb C_A^n$, has dimension ${{n-1}\choose 2}$.  A  basis for $\mathbb C^n_A$, which consists of  
three-cycles with equal costs between successive vertices, is 
\begin{equation}\label{def: A basis}  \{V_1 \stackrel{1}{\longrightarrow} V_j \stackrel{1}{\longrightarrow} V_k \stackrel{1}{\longrightarrow} V_1\}_{1<j<k\le n}. \end{equation}  \end{Thm} 

The
Eq.~\ref{def: A basis} cycles  
are anchored at one vertex, so the following offers a  more general  choice.  
\begin{Cor}\label{cor: general 3basis}  If $\mathcal{CB}^n_A$ has ${n-1}\choose2$ three-cycles where  each arc in a three-cycle has length 1 and  
each three-cycle has  one arc that is  
 not in any other  $\mathcal{CB}^n_A$ three-cycle, then $\mathcal{CB}^n_A$ is a $\mathbb C^n_A$   basis.  
 \end{Cor}

According to Thm.~\ref{thm: cn}, the  $\mathcal G^n_{A, cyclic} \in \mathbb C^n_A$ structure is governed by  
  three-cycles. To motivate their Eq.~\ref{def: A basis} form,    strong transitivity requires 
 $V_1  \stackrel{x}{\longrightarrow}V_j  \stackrel{y}{\longrightarrow} V_k = V_1  \stackrel{z}{\longrightarrow} V_k$, which defines the equation 
 $x+y-z=0$.  This equation  has the   normal vector $(1,1, -1)$, which,   
 when expressed in a path form,  
  is $V_1 \stackrel{1}{\longrightarrow} V_j \stackrel{1}{\longrightarrow} V_k \stackrel{1}{\longrightarrow} V_1$, or  an  Eq.~\ref{def: A basis} three-cycle. 
 It follows that the multiple of a three-cycle     measures how  
    this triplet's  $\mathcal G^n_A$ data portion  
 deviates from the triplet's cpi ``sameness."

  This discussion  leads to 
  the following central result where  
    Eq.~\ref{eq: splitting} asserts  that 
  the   Eq.~\ref{eq: first decomp}  goal     has been realized.   
  The theorem's      concluding statement      is crucial for what follows in this section.

\begin{Thm}\label{thm: gen decomp}   Space $\mathbb G^n_A$ is divided into a linear subspace $\mathbb{ST}^n_A$ and its orthogonal complement $\mathbb C^n_A$. For  $\mathcal G_A^n \in \mathbb G^n_A$, there are  unique  
 $\mathcal G^n_{A, cpi} \in \mathbb {ST}^n_A$  
  and   $\mathcal G^n_{A, cyclic} \in \mathbb C^n_A$ so that
\begin{equation}\label{eq: splitting} \mathcal G_A^n = \mathcal G^n_{A, cpi} + \mathcal G^n_{A, cyclic};\end{equation}  
 $\mathcal G^n_{A, cpi}$ and $\mathcal G^n_{A, cyclic}$ are, respectively, the orthogonal projections of $\mathcal G^n_A$ to $\mathbb{ST}_A^n$ and to $\mathbb C^n_A$.   The length of a closed path in  $\mathcal G_A^n$  equals    its $\mathcal G^n_{A, cyclic}$ length.\end{Thm}

The critical  last statement is a consequence of the  linear form of Eq.~\ref{eq: splitting}, which  
   requires the length of a  path in $\mathcal G_A^n$ to  equal the sum of its lengths in  $\mathcal G^n_{A, cpi}$ and in  $\mathcal G^n_{A, cyclic}.$ By design, the length of a closed path in  $\mathcal G^n_{A, cpi}$  is zero.  Namely,  
  $\mathcal G^n_{A, cpi}$     extracts   those   portions of $\mathcal G^n_A$ entries  that have
 {\em nothing substantive} to contribute  to closed path lengths.   
It now follows  
  that the path's lengths  in $\mathcal G_A^n$ and in $\mathcal G^n_{A, cyclic}$ must agree.  
In turn, this means that  
  {\em  all  relevant  closed path  information  for $\mathcal G^n_A$  
  is encoded in the three-cycles of $\mathcal G^n_{A, cyclic}$.}     Stated differently,  
  $\mathcal G^n_A\in \mathbb G^n_A$ has  
  an inherent three-cycle symmetry structure  
  displayed  by  
   $\mathcal G^n_{A, cyclic}$ but        camouflaged by     $\mathcal G^n_{A, cpi}$.   
 
 \begin{tikzpicture}[xscale=0.45, yscale=0.42]

\draw[->] (0, 0) -- (2.9, 0);
\draw[->]  (0, 0) -- (1.45, 2.78);
\draw[->] (3, 0) -- (1.55, 2.82); 
\node[below] at (0, 0) {$V_1$};
\node[below] at (3, 0) {$V_2$};
\node[above] at (1.5, 2.8) {$V_3$};
\node[below] at (1.5, 0) {8};
\node[left] at (.7, 1.5) {12};
\node[right] at (2.3, 1.5) {10};

\node[below] at (1.5, -1) {{\bf a.} $\mathcal G^3_A$};

\node at (4.5, 1.5) {=};

%%%%%%%%%%%%%

\draw[->] (6, 0) -- (8.9, 0);
\draw[->]  (6, 0) -- (7.45, 2.78);
\draw[->] (9, 0) -- (7.55, 2.82); 
\node[below] at (6, 0) {$V_1$};
\node[below] at (9, 0) {$V_2$};
\node[above] at (7.5, 2.8) {$V_3$};
\node[below] at (7.5, 0) {6};
\node[left] at (6.7, 1.5) {14};
\node[right] at (8.3, 1.5) {8};

\node[below] at (7.5, -1) {{\bf b.} $\mathcal G^3_{A, cpi}$};

\node at (10.5, 1.5) {+};

%%%%%%%%%%%%%%%%%

\draw[->] (12, 0) -- (14.9, 0);
\draw[<-]  (12.1, 0.1) -- (13.45, 2.78);
\draw[->] (15, 0) -- (13.55, 2.82); 
\node[below] at (12, 0) {$V_1$};
\node[below] at (15, 0) {$V_2$};
\node[above] at (13.5, 2.8) {$V_3$};
\node[below] at (13.5, 0) {2};
\node[left] at (12.7, 1.5) {2};
\node[right] at (14.3, 1.5) {2};

\node[below] at (13.5, -1) {{\bf c.} $\mathcal G^3_{A, cyclic}$};

\node at (7.5, -2.8) {{\bf Figure 2.} Interpreting $\mathcal G^n_{A, cyclic}$};

\end{tikzpicture}

 To expand on the comment that  the $\mathcal G^n_{A, cpi}$ entries contribute {\em nothing substantive} about  
 closed path lengths,  
  notice that   computing  the  
     $V_1 \stackrel{8}{\longrightarrow}   V_2  \stackrel{10}{\longrightarrow}   V_3  \stackrel{-12}{\longrightarrow}   V_1$  length of 6  in Fig.~2a 
  involves a subtraction cancelation. To appreciate this structure,  
 let the optimal  (but unknown)  cancelled values  be $u$, $v$, and $w$ from, respectively, arcs $\widehat{V_1V_2}$, $\widehat{V_2V_3}$, and $\widehat{V_3V_1}$.     
  That is, $(8-u)+(10-v)+(-12-w)=6$ where the cancelled values define the  equation     $u+v+w=0$, which corresponds to  a zero-length closed path.  
  For $n>3$, this cancellation  applies to all triplets, so  these extracted   
   values   define a  $\mathbb{ST}_A^n$ graph.  The optimal  
   choice of removed terms comes from the $\mathbb{ST}_A^n$ graph that most closely resembles $\mathcal G^n_A$, which is its orthogonal projection $\mathcal G^n_{A, cpi}$ (Thm.~\ref{thm: gen decomp}).      Indeed, with Fig.~2, the $\mathcal G^3_{A, cpi}$ component extracts  
    $u+v+w= 6+8-14=0$.  What remains are  portions of arc entries that, without  further  need of modification, are   relevant   for computing path lengths.   
     These terms define     the $\mathcal G^3_{A, cyclic}$ graph with its  $V_1  \stackrel{2}{\longrightarrow}   V_2  \stackrel{2}{\longrightarrow}   V_3  \stackrel{2}{\longrightarrow}  V_1$ closed cycle  that  
   directly provides the path length of 6.   
   Not only do the superfluous  $\mathcal G^n_{A, cpi}$ terms complicate computations, but, as discussed  below,      they can sidetrack  optimization approaches such as the greedy algorithm.

The same behavior  holds in general; e.g., the zero length of each triplet in  Fig.~1b   
  identifies the optimal subtraction cancelations   
   for  computing  $\mathcal G^5_A$ path lengths.  As  the cyclic $\mathcal G^n_{A, cyclic}$  captures the germane portions for path-length considerations, when searching for the longest or shortest $\mathcal G^n_A$  paths, ignore    
  $\mathcal G^n_A$ and   analyze  only   the simpler    $\mathcal G^n_{A, cyclic}$.   
  (Some $\mathcal G^n_{A, cyclic}$ arcs  
    belong to several cycles; e.g.,   the $V_2  \stackrel{2}{\longrightarrow}  V_1$ leg in Fig.~1c is  the sum of this arc's length in two cycles.)

 A slight modification of Thm.~\ref{thm: gen decomp}    describes the length of  any connected path.   
  
 \begin{Cor}\label{cor: general path} The length of a path in $\mathcal G^n_A$ that connects $V_j$ with $V_k$ is the length of this path in $\mathcal G^n_{A, cyclic}$ plus the length of the $V_j\to V_k$ arc in $\mathcal G^n_{A, cpi}$.\end{Cor}

  To illustrate, the length of  12 for path  $V_1 \stackrel{14}{\longrightarrow}  V_4 \stackrel{-6}{\longrightarrow}  V_3 \stackrel{1}{\longrightarrow}   V_5 \stackrel{3}{\longrightarrow}  V_4$ in $\mathcal G^5_A$ (Fig.~1a), which can meet vertices multiple times, equals the  easier computed  length of 1 for    $V_1 \stackrel{3}{\longrightarrow}  V_4 \stackrel{0}{\longrightarrow}  V_3 \stackrel{-2}{\longrightarrow}   V_5 \stackrel{0}{\longrightarrow}  V_4$    in $\mathcal G^5_{A, cyclic}$ plus 11 from the $V_1\stackrel{11}{\longrightarrow}  V_4$ arc length in $\mathcal G^5_{A, cpi}$.  The $3-2$ subtraction manifests the algebraic  arrangement of the triplets.\smallskip
 
 {\em Proof:} The length of a path in $\mathcal G^n_A$ equals the sum of its lengths in $\mathcal G^n_{A, cpi}$ and $\mathcal G^n_{A, cyclic}$.  According to Cor.~\ref{cor: gen path length}, the length of a connected path in $\mathcal G^n_{A, cpi}$ starting at $V_j$ and ending at $V_k$ equals the length of the   arc $V_j\to V_k$ connecting the endpoints.  This completes the proof. $\square$\smallskip

In general,   
it  is easier to
extract  $\mathcal G^n_A$ path properties   
from  $\mathcal G^n_{A, cyclic}$ than from $\mathcal G^n_A$; e.g.,   even the flawed   
 greedy algorithm (GA)   shows that $V_1 \stackrel{3}{\longrightarrow}  V_4 \stackrel{3}{\longrightarrow} V_2 \stackrel{2}{\longrightarrow}  V_5 \stackrel{2}{\longrightarrow} V_3 \stackrel{1}{\longrightarrow}  V_1$ of length 11 is  the longest $\mathcal G^5_{A, cyclic}$  (Fig.~1c)  Hamiltonian path.  According to Eq.~\ref{def: dij=-dji},   its reversal (length -11) is the shortest.   This solves the  $\mathcal G^5_A$ TSP  problem because, as  Thm.~\ref{thm: gen decomp}   asserts,      these two $\mathcal G^5_{A, cyclic}$  paths identify, respectively,       the longest and shortest $\mathcal G^5_{A}$ Hamiltonian paths and  their lengths.  But GA\footnote{Because its limitations and failings are  well known, the   greedy algorithm (GA) is used   to illustrate advantages of the decomposition. Other GA difficulties, caused by the algebra of cycles, are indicated with Fig.~12.} fails with $\mathcal G^5_A$ primarily because  the $\mathcal G^5_{A, cpi}$    values, which  cancel when calculating  lengths,  divert the algorithm.   To see this, using the GA to search for the longest Hamiltonian path of $\mathcal G^5_A$ yields the incorrect $V_1 \stackrel{14}{\longrightarrow}   V_4 \stackrel{-3}{\longrightarrow} V_5 \stackrel{-1}{\longrightarrow} V_3 \stackrel{-3}{\longrightarrow} V_2 \stackrel{1}{\longrightarrow}  V_1$.    
   Applying GA to $\mathcal G^5_{A, cpi}$ generates the same path, which  underscores    the fact  
    that  $\mathcal G^5_{A, cpi}$ is the source of this problem.

 For $\mathbb G^n_A$, the tasks of finding the longest and shortest $\mathcal G^n_{A, cyclic}$  Hamiltonian paths coincide.  This is because the reversal of one is the other.
 
 \begin{Cor} If the length of a  path in $\mathcal G^n_{A, cyclic}$ is $x$, then the length of its reversal is $-x$.  \end{Cor}

 \subsection{Decomposition}\label{sect: compute}  To compute  
  $\mathcal G^n_{A, cpi}$ and $\mathcal G^n_{A, cyclic}$,   recall  that $\mathcal G^n_{A, cpi}\in \mathbb{ST}^n_A$   is the    orthogonal projection of $\mathcal G^n_A$.   
   Terms from this projection are described next; for added details  see \cite{tandd}.    

The projection is a  
linear algebra exercise.  To connect graphs with vectors, let $\mathbf d_A^n\in \mathbb R_A^{n\choose2}$ be

 \begin{equation}\label{def: dn}\mathbf d_A^n=(d_{1, 2}, d_{1, 3}, \dots, d_{1, n}; d_{2, 3}, \dots, d_{2, n}; d_{3, 4}, \dots ; d_{n-1, n}), \textrm{ where } d_{i, j} = -d_{j, i};\end{equation}

\noindent  the semicolons designate where the first subscript changes.  To identify $\mathbb G^n_A$ with     $\mathbb R_A^{n\choose2}$,   
 let $d_{i, j}$ be the arc length $V_i  \stackrel{d_{i, j}}{\longrightarrow}  V_j.$  As  $V_i \stackrel{d_{i, j}}{\longrightarrow}  V_j$  equals    $V_j \stackrel{-d_{i, j}}{\longrightarrow}  V_i$, it follows that $d_{j, i}=-d_{i, j}$.     Let   $\mathbb{ST}_A^n \subset \mathbb G^n_A$    also denote the   $(n-1)$-dimensional (strongly transitive) 
subspace of $\mathbb R_A^{n\choose2}$ where  
 each triplet $\{i, j, k\}$  satisfies $d_{i, j} + d_{j, k} = d_{i, k}$.  With these identifications, structures  of $\mathbb R^{n\choose2}_A$  and  $\mathbb G^n_A$ can  be described interchangeably.    
 
 \begin{Def}\label{def: S} For  vertex  $V_j$  
 of $\mathcal G_A^n \in \mathbb G^n_A$,   let  
 $\mathcal S_A(V_j)$ be  $\frac1n$ times the sum of the arc lengths  leaving  
 vertex $V_j$, $j=1, \dots, n.$   
 \end{Def}

 \begin{Thm}\label{thm: computing cpi} \cite{tandd} For $\mathcal G_A^n\in \mathbb G^n_A$,   the $\mathcal G^n_{A, cpi}$  path length from $V_i$ to $V_j$ is  
  $d_{i, j}=\mathcal S_A(V_i) - \mathcal S_A(V_j)$, $ i, j \in \{1, \dots, n\}$.  Each of $\mathcal G_A^n$ and $\mathcal G_{A, cpi}^n$ satisfy  
   $\sum_{j=1}^n \mathcal S_A(V_j)=0$.   Graph  $\mathcal G^n_{A, cyclic}$ is given by  $\mathcal G^n_{A, cyclic} =\mathcal G^n_A-\mathcal G^n_{A, cpi}$.
     All  vertices of  $\mathcal G^n_{A, cyclic}$  
    satisfy the stronger  $ \mathcal S_A(V_j)=0. $   Conversely, if all vertices of $\mathcal G^n_A\in \mathbb G^n_A$ satisfy $ \mathcal S_A(V_j)=0, $ then $\mathcal G^n_A\in \mathbb C^n_A$ 
  \end{Thm}\smallskip
  
The concluding  statement follows from Thm.~\ref{thm: computing cpi}'s first sentence.  This is  because $ \mathcal S_A(V_j)=0 $  for all vertices requires all $d_{i, j}$  legs of $\mathcal G^n_{A, cpi}$ to equal zero.  As $\mathcal G^n_{A, cpi}=0$,  $\mathcal G^n_A$  equals its $\mathcal G^n_{A, cyclic}$ component.

With $\mathcal G^6_{A}$ (Fig.~3a), the $S_A(V_j)$ values  (called  `Borda Values' in \cite{tandd})  are   
  $ \mathcal S_A(V_1) = \frac16(-1-5 -3+4 -1) =-1,  \mathcal S_A(V_2) = 0, \mathcal S_A(V_3) = 2, \mathcal S_A(V_4)= -2, \mathcal S_A(V_5)=-3, \mathcal S_A(V_6) = 4.$  
  Thus (Thm.~\ref{thm: computing cpi}),  
  the     Fig.~3b values for $\mathcal G^6_{A, cpi}$  are 
 $d_{1, 2} =  \mathcal S_A(V_1) -  \mathcal S_A(V_2) = -1,  d_{1, 3} = 3, d_{1, 4} = 1, d_{1, 5} = 2, d_{1, 6}=-5, d_{2, 3} = -2, d_{2, 4} = 2,   d_{2, 5} = 3,  d_{2, 6} = -4, d_{3, 4} = 4,  d_{3, 5}=5, d_{3, 6}= -2, d_{4, 5} = 1, d_{4, 6}=-6, d_{5, 6}= -7.$ 
 Graph $\mathcal G^6_{A, cyclic}$   
 follows from the equality   $\mathcal G^6_{A, cyclic} = \mathcal G^6_A - \mathcal G^6_{A, cpi}$; this 
   defines  
 Fig.~3c.  
 
   \begin{tikzpicture}[xscale=0.5, yscale=0.52]

\draw[<-] (1.18, 0) -- (3.9, 0); %ab
\node[below] at (1, 0) {$V_1$};
\node[below] at (2.5, 0) {1};
\draw[<-] (1.15, .06) -- (4.9, 2);   %ac
\node at (1.9, .5) {\small 5};
\draw[<-] (1.15, .15) -- (3.9, 3.9); %ad
\node at (2.1, 1.5)  {\small 3};
\draw[->] (1, .2) -- (1, 3.9);  %a, e)
\node at  (1, 3) {\small 4}; 
\node[left] at (0, 2) {$V_6$}; 
\draw[<-] (.9, .1) -- (.1, 2);  %af
\node[left] at (.5, 1) {1};
\draw[<-] (4.1, .1) -- (4.9, 1.9); %bc
\node[right] at  (4.5, 1) {7}; 
\node[below] at (4, 0) {$V_2$};
\draw[->] (4, .1) -- (4, 3.85);  %bd
\node at (4, 1) {\small 7};
\draw[->]  (3.9, .1) -- (1.1, 3.9); %be
\node at (2.1, 2.5) {\small 3};
\draw[<-] (3.82, .07) -- (.23, 1.95);  %bf
\node at (.56, 1.7) {\small 4};
\draw[->] (4.9, 2.1) -- (4.1, 3.9);  %cd
\node[right] at (4.5, 3) {2}; 
\node[right] at (5, 2){$V_3$};
\draw[->] (4.9, 2.05) -- (1.24, 3.86); %ce
\node at (4.3, 2.35) {\small 3};
\draw[<-] (4.8, 2) -- (.2, 2); %cf
\node at (3.3, 2) {\small 5};
\draw[->] (3.8, 4) -- (1.2, 4); %de
\node[above] at (2.5, 4) {1};
\node[above] at (4, 4) {$V_4$}; 
\draw[<-] (3.75, 3.85) -- (.2, 2.1);  %df
\node at (3.05, 3.4)   {\small 7};
\node[above] at (1, 4) {$V_5$}; 
\draw[<-]  (.83, 3.9) -- (.1, 2.1); %ef
\node[left] at (.5, 3) {7};

\node[below] at (2, -1) {{\bf a.} $\mathcal G_A^6$};

%%%%%%%%%%%%%%%%%

\draw[<-] (11.18, 0) -- (13.9, 0); %ab
\node[below] at (11, 0) {$V_1$};
\node[below] at (12.5, 0) {1};
\draw[<-] (11.15, .06) -- (14.9, 2);   %ac
\node at (11.9, .5) {\small 3};
\draw[->] (11.15, .15) -- (13.9, 3.9); %ad
\node at (12.1, 1.5)  {\small 1};
\draw[->] (11, .2) -- (11, 3.9);  %a, e)
\node at  (11, 3) {\small 2}; 
\node[left] at (10, 2) {$V_6$}; 
\draw[<-] (10.9, .1) -- (10.1, 2);  %af
\node[left] at (10.5, 1) {5};
\draw[<-] (14.1, .1) -- (14.9, 1.9); %bc
\node[right] at  (14.5, 1) {2}; 
\node[below] at (14, 0) {$V_2$};
\draw[->] (14, .1) -- (14, 3.85);  %bd
\node at (14, 1) {\small 2};
\draw[->]  (13.9, .1) -- (11.1, 3.9); %be
\node at (12.1, 2.5) {\small 3};
\draw[<-] (13.82, .07) -- (10.23, 1.95);  %bf
\node at (10.56, 1.7) {\small 4};
\draw[->] (14.9, 2.1) -- (14.1, 3.9);  %cd
\node[right] at (14.5, 3) {4}; 
\node[right] at (15, 2){$V_3$};
\draw[->] (14.9, 2.05) -- (11.24, 3.86); %ce
\node at (14.3, 2.35) {\small 5};
\draw[<-] (14.8, 2) -- (10.2, 2); %cf
\node at (13.3, 2) {\small 2};
\draw[->] (13.8, 4) -- (11.2, 4); %de
\node[above] at (12.5, 4) {1};
\node[above] at (14, 4) {$V_4$}; 
\draw[<-] (13.75, 3.85) -- (10.2, 2.1);  %df
\node at (13.05, 3.45)   {\small 6};
\node[above] at (11, 4) {$V_5$}; 
\draw[<-]  (10.83, 3.9) -- (10.11, 2.1); %ef
\node[left] at (10.5, 3) {7};

\node at (7.5, 2) {=};

\node[below] at (12, -1) {{\bf b.} $\mathcal G_{A, cpi}^6$};

%%%%%%%%%%%%%

\draw[dashed] (21.18, 0) -- (23.9, 0); %ab
\node[below] at (21, 0) {$V_1$};
\node[below] at (22.5, 0) {0};
\draw[<-] (21.15, .06) -- (24.9, 2);   %ac
\node at (21.9, .5) {\small 2};
\draw[<-] (21.15, .15) -- (23.9, 3.9); %ad
\node at (22.1, 1.5)  {\small 4};
\draw[->] (21, .2) -- (21, 3.9);  %a, e)
\node at  (21, 3) {\small 2}; 
\node[left] at (20, 2) {$V_6$}; 
\draw[->] (20.9, .1) -- (20.1, 2);  %af
\node[left] at (20.5, 1) {4};
\draw[<-] (24.1, .1) -- (24.9, 1.9); %bc
\node[right] at  (24.5, 1) {5}; 
\node[below] at (24, 0) {$V_2$};
\draw[->] (24, .1) -- (24, 3.85);  %bd
\node at (24, 1) {\small 5};
\draw[dashed]  (23.9, .1) -- (21.1, 3.9); %be
\node at (22.1, 2.5) {\small 0};
\draw[dashed] (23.82, .07) -- (20.23, 1.95);  %bf
\node at (20.56, 1.7) {\small 0};
\draw[<-] (24.9, 2.1) -- (24.1, 3.9);  %cd
\node[right] at (24.5, 3) {2}; 
\node[right] at (25, 2){$V_3$};
\draw[<-] (24.9, 2.05) -- (21.24, 3.86); %ce
\node at (24.3, 2.35) {\small 2};
\draw[<-] (24.8, 2) -- (20.2, 2); %cf
\node at (23.3, 2) {\small 3};
\draw[dashed]  (23.8, 4) -- (21.2, 4); %de
\node[above] at (22.5, 4) {0};
\node[above] at (24, 4) {$V_4$}; 
\draw[<-] (23.75, 3.85) -- (20.2, 2.1);  %df
\node at (23.05, 3.4)   {\small 1};
\node[above] at (21, 4) {$V_5$}; 
\draw[dashed]  (20.83, 3.9) -- (20.11, 2.1); %ef
\node[left] at (20.5, 3) {0};

\node at (17.5, 2) {+};

\node[below] at (22, -1) {{\bf c.} $\mathcal G_{A, cyclic}^6$};

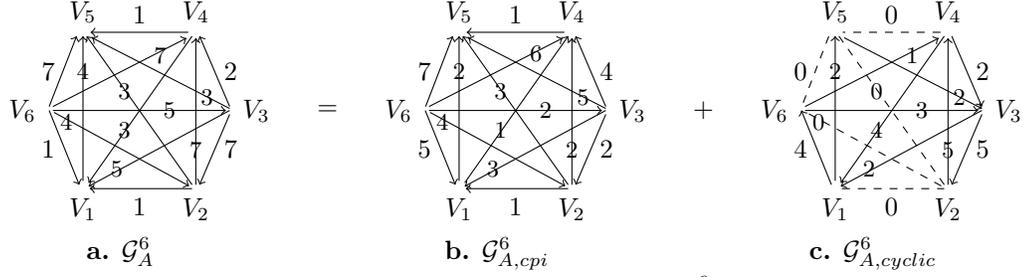
 \node[below] at (12.5, -2) {{\bf Figure 3.} Decomposition of a $\mathcal G_A^6$};

\end{tikzpicture}

 Notice how the redundant   $\mathcal G^6_{A, cpi}$ (Fig.~3b)   dominates the $\mathcal G_A^6$  structure 
 even though    the     simpler   
   $\mathcal G^6_{A, cyclic}$ (Fig.~3c)  
is what   determines all of   $\mathcal G_A^6$'s closed path properties.   This must be expected because, according to the orthogonal projection construction,  $\mathcal G^n_{A, cpi}$ is the $\mathbb {ST}_A^n$ graph that most closely resembles $\mathcal G^n_A$.     As true with Fig.~2,  a feature of  this  $\mathcal G^6_{A, cpi}$ and $\mathcal G^6_A$ similarity is that  $\mathcal G^6_{A, cpi}$   collects  terms involved in subtraction/cancellations when computing $\mathcal G^6_A$  path lengths; the remaining $\mathcal G^6_{A, cyclic}$ entries are the relevant  portions for determining path properties.

 \begin{tikzpicture}[xscale=0.5, yscale=0.53]

\draw[->] (1.18, 0) -- (3.9, 0); %ab
\node[below] at (1, 0) {$V_1$};
\node[below] at (2.5, 0) {3};
\draw[->] (1.15, .06) -- (4.9, 2);   %ac
\node at (1.9, .5) {\small 4};
\draw[->] (1.15, .15) -- (3.9, 3.9); %ad
\node at (2.1, 1.5)  {\small 15};
\draw[->] (1, .2) -- (1, 3.9);  %a, e)
\node at  (1, 3) {\small 21}; 
\node[left] at (0, 2) {$V_6$}; 
\draw[->] (.9, .1) -- (.1, 2);  %af
\node[left] at (.5, 1) {27};
\draw[->] (4.1, .1) -- (4.9, 1.9); %bc
\node[right] at  (4.5, 1) {3}; 
\node[below] at (4, 0) {$V_2$};
\draw[->] (4, .1) -- (4, 3.85);  %bd
\node at (4, 1) {\small 10};
\draw[->]  (3.9, .1) -- (1.1, 3.9); %be
\node at (2.1, 2.5) {\small 17};
\draw[->] (3.82, .07) -- (.23, 1.95);  %bf
\node at (.56, 1.7) {\small 25};
\draw[->] (4.9, 2.1) -- (4.1, 3.9);  %cd
\node[right] at (4.5, 3) {13}; 
\node[right] at (5, 2){$V_3$};
\draw[->] (4.9, 2.05) -- (1.24, 3.86); %ce
\node at (4.3, 2.35) {\small 16};
\draw[->] (4.8, 2) -- (.2, 2); %cf
\node at (3.3, 2) {\small 24};
\draw[->] (3.8, 4) -- (1.2, 4); %de
\node[above] at (2.5, 4) {5};
\node[above] at (4, 4) {$V_4$}; 
\draw[->] (3.75, 3.85) -- (.2, 2.1);  %df
\node at (3.05, 3.4)   {\small 13};
\node[above] at (1, 4) {$V_5$}; 
\draw[->]  (.83, 3.9) -- (.1, 2.1); %ef
\node[left] at (.5, 3) {9};

\node[below] at (2, -1) {{\bf a.} $\mathcal G_A^6$};

%%%%%%%%%%%%%%%%%%%%

\draw[->] (11.18, 0) -- (13.9, 0); %ab
\node[below] at (11, 0) {$V_1$};
\node[below] at (12.5, 0) {3};
\draw[->] (11.15, .06) -- (14.9, 2);   %ac
\node at (11.9, .5) {\small 4};
\draw[->] (11.15, .15) -- (13.9, 3.9); %ad
\node at (12.1, 1.5)  {\small 15};
\draw[->] (11, .2) -- (11, 3.9);  %a, e)
\node at  (11, 3) {\small 20}; 
\node[left] at (10, 2) {$V_6$}; 
\draw[->] (10.9, .1) -- (10.1, 2);  %af
\node[left] at (10.5, 1) {28};
\draw[->] (14.1, .1) -- (14.9, 1.9); %bc
\node[right] at  (14.5, 1) {1}; 
\node[below] at (14, 0) {$V_2$};
\draw[->] (14, .1) -- (14, 3.85);  %bd
\node at (14, 1) {\small 12};
\draw[->]  (13.9, .1) -- (11.1, 3.9); %be
\node at (12.1, 2.5) {\small 17};
\draw[->] (13.82, .07) -- (10.23, 1.95);  %bf
\node at (10.56, 1.7) {\small 25};
\draw[->] (14.9, 2.1) -- (14.1, 3.9);  %cd
\node[right] at (14.5, 3) {11}; 
\node[right] at (15, 2){$V_3$};
\draw[->] (14.9, 2.05) -- (11.24, 3.86); %ce
\node at (14.3, 2.35) {\small 16};
\draw[->] (14.8, 2) -- (10.2, 2); %cf
\node at (13.3, 2) {\small 24};
\draw[->] (13.8, 4) -- (11.2, 4); %de
\node[above] at (12.5, 4) {5};
\node[above] at (14, 4) {$V_4$}; 
\draw[->] (13.75, 3.85) -- (10.2, 2.1);  %df
\node at (13.05, 3.4)   {\small 13};
\node[above] at (11, 4) {$V_5$}; 
\draw[->]  (10.83, 3.9) -- (10.1, 2.1); %ef
\node[left] at (10.5, 3) {8};

\node[below] at (12, -1) {{\bf b.} $\mathcal G_{A, cpi}^6$};

%%%%%%%%%%%%%%%%%

\draw[dashed] (21.18, 0) -- (23.9, 0); %ab
\node[below] at (21, 0) {$V_1$};
\node[below] at (22.5, 0) {0};
\draw[dashed] (21.15, .06) -- (24.9, 2);   %ac
\node at (21.9, .5) {\small 0};
\draw[dashed] (21.15, .15) -- (23.9, 3.9); %ad
\node at (22.1, 1.5)  {\small 0};
\draw[thick, ->] (21, .2) -- (21, 3.9);  %a, e)
\node at  (20.9, 3) {\small 1}; 
\node[left] at (20, 2) {$V_6$}; 
\draw[thick, <-] (20.9, .1) -- (20.1, 2);  %af
\node[left] at (20.5, 1) {1};
\draw[thick, ->] (24.1, .1) -- (24.9, 1.9); %bc
\node[right] at  (24.5, 1) {2}; 
\node[below] at (24, 0) {$V_2$};
\draw[thick, <-] (24, .1) -- (24, 3.85);  %bd
\node at (24, 1) {\small 2};
\draw[dashed]  (23.9, .1) -- (21.1, 3.9); %be
\node at (22.1, 2.5) {\small 0};
\draw[dashed] (23.82, .07) -- (20.23, 1.95);  %bf
\node at (20.56, 1.7) {\small 0};
\draw[thick, ->] (24.9, 2.1) -- (24.1, 3.9);  %cd
\node[right] at (24.5, 3) {2}; 
\node[right] at (25, 2){$V_3$};
\draw[dashed] (24.9, 2.05) -- (21.24, 3.86); %ce
\node at (24.3, 2.35) {\small 0};
\draw[dashed] (24.8, 2) -- (20.2, 2); %cf
\node at (23.3, 2) {\small 0};
\draw[dashed] (23.8, 4) -- (21.2, 4); %de
\node[above] at (22.5, 4) {0};
\node[above] at (24, 4) {$V_4$}; 
\draw[dashed] (23.75, 3.85) -- (20.2, 2.1);  %df
\node at (23.05, 3.45)   {\small 0};
\node[above] at (21, 4) {$V_5$}; 
\draw[thick,->]  (20.83, 3.9) -- (20.11, 2.1); %ef
\node[left] at (20.5, 3) {1};

 \node at (7.5, 2) {$=$};
 \node at (17.5, 2) {$+$};

\node[below] at (23, -1) {{\bf c.} $\mathcal G_{A, cyclic}^6$};

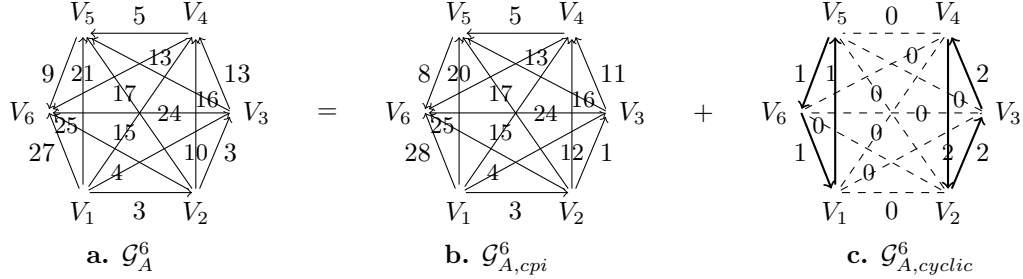
 \node[below] at (12.5, -2) {{\bf Figure 4.} Advantages of  $\mathcal G_{A, cyclic}^n$};

\end{tikzpicture}

    Figure 4 illustrates  
    how seriously the $\mathcal G^n_{A, cpi}$ terms can  cloud a path analysis.      While path properties are determined by the extremely simple  
      $\mathcal G^6_{A, cyclic}$  (Fig.~4c), this clarity  is not obvious from  $\mathcal G^6_A$ (Fig.~4a).  The reason is that    $\mathcal G^6_A$  more  closely resembles the associated $\mathcal G^6_{A, cpi}$.
As developed next,   this is a general phenomenon.

\begin{Def}\label{def: A equivalence} Two graphs $\mathcal G^n_{A, 1}, \mathcal G^n_{A, 2} \in \mathbb G^n_A$   are ``closed path  equivalent" iff $\mathcal G^n_{A, 1, cpi} = \mathcal G^n_{A, 2, cpi}$.  \end{Def} 

\begin{Cor} \label{cor: cpe} The ``closed path equivalent" relationship is an equivalence relation.  
 Two graphs are equivalent iff their difference is a graph in $\mathbb{ST}_A^n$.  Thus, an equivalence class of this relationship is the sum of a  $\mathcal G^n_{A, cyclic}$ and the $(n-1)$-dimensional linear subspace $\mathbb {ST}_A^n$. \end{Cor}

According to  Cor.~\ref{cor: cpe},   multiple (actually, most)  choices     of a $\mathcal G^n_{A, cpi}$ from the vast offerings of the $(n-1)$-dimensional $\mathbb{ST}_A^n$   dictate the form of  $\mathcal G^n_A$ and  obscure the relevant $\mathcal G^n_{A, cyclic}$.

\subsection{Structure of triplets}  
Analyzing   closed path properties of a $\mathcal G^n_A$ involves  
  the algebraic structure of the $\mathcal G^n_{A, cyclic}$ three-cycles.       The following  shows how    to  identify    $\mathcal G^n_{A, cyclic}$'s   
  three-cycles.  

\begin{Thm}\label{thm: finding cycles}  For $\mathcal G^n_{A, cyclic}$,  only one three-cycle of a Cor.~\ref{cor: general 3basis} basis has a $\widehat{V_jV_k}$ arc.  The cycle's  multiple is    
   the $V_j  \stackrel{d_{j, k}}{\longrightarrow} V_k$ weight in $\mathcal G^n_{A, cyclic}$.\end{Thm}

\noindent{\em Proof:}  
Arc $V_j  \stackrel{d_{j, k}}{\longrightarrow} V_k$ in   $\mathcal G^n_{A, cyclic}$  appears only in   $V_s \stackrel{x}{\longrightarrow} V_j \stackrel{x}{\longrightarrow} V_k \stackrel{x}{\longrightarrow} V_s$  of the $\mathcal {CB}^n_A$  
basis. As all  weights in a three-cycle agree,  this is the cycle's multiple.  
    $\square$\smallskip

To illustrate,   Fig.~1c has the three three-cycles $V_1 \stackrel{3}{\longrightarrow} V_4 \stackrel{3}{\longrightarrow} V_2 \stackrel{3}{\longrightarrow}  V_1, \, V_1 \stackrel{1}{\longrightarrow}  V_2 \stackrel{1}{\longrightarrow}  V_3 \stackrel{1}{\longrightarrow}  V_1$ and $V_2 \stackrel{2}{\longrightarrow} V_5 \stackrel{2}{\longrightarrow} V_3 \stackrel{2}{\longrightarrow} V_2$. 
 Figure~3c has the four three-cycles  
 $V_2 \stackrel{5}{\longrightarrow} V_4 \stackrel{5}{\longrightarrow} V_3 \stackrel{5}{\longrightarrow} V_2,  V_1 \stackrel{4}{\longrightarrow} V_6 \stackrel{4}{\longrightarrow} V_4 \stackrel{4}{\longrightarrow} V_1,   V_3 \stackrel{3}{\longrightarrow} V_4 \stackrel{3}{\longrightarrow} V_6 \stackrel{3}{\longrightarrow} V_3,$ and  $V_1 \stackrel{2}{\longrightarrow} V_5 \stackrel{2}{\longrightarrow} V_3 \stackrel{2}{\longrightarrow} V_1.$   When finding Hamiltonian paths,   to avoid prematurely  returning to a vertex,    at most two arcs of a three-cycle can be used.  With this caveat, the GA can succeed with  
   $\mathcal G^n_{A, cyclic}$  settings where it would fail   for $\mathcal G^n_A$.  For instance, the GA delivers  
   the longest  Fig.~3c     Hamiltonian  path      
 $ V_3 \stackrel{5}{\longrightarrow} V_2 \stackrel{5}{\longrightarrow} V_4 \stackrel{4}{\longrightarrow} V_1 \stackrel{4}{\longrightarrow} V_6 \stackrel{0}{\longrightarrow} V_5 \stackrel{2}{\longrightarrow} V_3$  
with   length $20,$   
 which   is its $\mathcal G^6_A$ (Fig.~3a) path length (Thm.~\ref{thm: gen decomp}).  For    $\mathcal G^n_A\in \mathbb G^n_A$, its shortest Hamiltonian route  
    reverses  the longest.

 A vertex that is a source or sink imposes an obstacle in finding optimal paths.  For instance,   if $V_j$ is a sink for negative cost directions, as is $V_6$ in Fig.~3a,  
 all ways to leave this vertex  require using a
   positive cost direction.  
    Fortunately,  general properties of  $\mathcal G^n_{A, cyclic} \in \mathbb C^n_A$ can be obtained via  the    $\mathbb C^n_A$ basis (Thm.~\ref{thm: cn}, Cor.~\ref{cor: general 3basis});       sample   conclusions    are in Thm.~\ref{thm: 4isgood}.  The first assertion identifies two well behaved settings.       The second comment asserts that although sources and sinks are not unusual in $\mathcal G^n_A$, which   cause subtraction cancellations, they never arise  in $\mathcal G^n_{A, cyclic}$.  According to the theorem's  last comment,  
    sources and sinks  appear  in $\mathcal G^n_A$ only because  $\mathcal G^n_{A, cpi}$ almost always has  them.  Thus if  $\mathcal G^n_A$ has a source and/or sink, expect that  its structure is dominated by  $\mathcal G^n_{A, cpi}$.     

\begin{Thm}\label{thm: 4isgood} For $\mathcal G^n_{A, cyclic} \in \mathbb C^n_A$, $n=4, 5$,  its longest and shortest Hamiltonian paths have, respectively,  all non-negative arc costs and non-positive arc costs.   For  $n\ge 4,$ while $\mathcal G^n_A$ can have sinks and/or sources, this  is impossible for a $\mathcal G^n_{A, cyclic}$.  In contrast, if all positive cost directions of $\mathcal G^n_{A, cpi}$ are non-zero, then both the positive  
and negative cost directions have a sink and a source.\end{Thm}

\begin{tikzpicture}[xscale=0.5, yscale=0.52]

\draw[->] (1.18, 0) -- (3.9, 0); %ab
\node[below] at (1, 0) {$V_1$};
\node[below] at (2.5, 0) {5};
\draw[->] (1.15, .06) -- (4.85, 1.95);   %ac
\node at (1.9, .5) {\small 5};
\draw[->] (1.15, .15) -- (3.8, 3.8); %ad
\node at (2.1, 1.5)  {\small 16};
\draw[->] (1, .2) -- (1, 3.9);  %a, e)
\node at  (1, 3) {\small 11}; 
\node[left] at (0, 2) {$V_6$}; 
 
\draw[<-] (4.1, .1) -- (4.9, 1.9); %bc
\node[right] at  (4.5, 1) {5}; 
\node[below] at (4, 0) {$V_2$};
\draw[->] (4, .1) -- (4, 3.85);  %bd
\node at (4, 1) {\small 17};
\draw[<-]  (3.9, .1) -- (1.1, 3.9); %be
\node at (2.1, 2.5) {\small 1};
\draw[<-] (3.82, .07) -- (.23, 1.95);  %bf
\node at (.56, 1.7) {\small 5};
 
\node[right] at (5, 2){$V_3$};
\draw[<-] (4.9, 2.05) -- (1.24, 3.86); %ce
\node at (4.3, 2.35) {\small 5};
 \node[above] at (4, 4) {$V_4$}; 
\draw[<-] (3.75, 3.85) -- (.2, 2.1);  %df
\node at (3.05, 3.4)   {\small 8};
\node[above] at (1, 4) {$V_5$}; 
\draw[<-]  (.83, 3.9) -- (.1, 2.1); %ef
\node[left] at (.5, 3) {12};

\node[below] at (2, -1) {{\bf a.} $\mathcal G_A^6$};

%%%%%%%%%%%%%%%%%

\draw[->] (11.18, 0) -- (13.9, 0); %ab
\node[below] at (11, 0) {$V_1$};
\node[below] at (12.5, 0) { 6};
\draw[->] (11.15, .06) -- (14.9, 2);   %ac
\node at (11.9, .5) {\small 7};
\draw[->] (11.15, .15) -- (13.9, 3.9); %ad
\node at (12.1, 1.5)  {\small 13};
\draw[->] (11, .2) -- (11, 3.9);  %a, e)
\node at  (11, 3) {\small 9}; 
\node[left] at (10, 2) {$V_6$}; 
\draw[->] (10.9, .1) -- (10.1, 2);  %af
\node[left] at (10.5, 1) {2};
\draw[->] (14.1, .1) -- (14.9, 1.9); %bc
\node[right] at  (14.5, 1) {1}; 
\node[below] at (14, 0) {$V_2$};
\draw[->] (14, .1) -- (14, 3.85);  %bd
\node at (14, 1) {\small 7};
\draw[->]  (13.9, .1) -- (11.1, 3.9); %be
\node at (12.1, 2.5) {\small 3};
\draw[<-] (13.82, .07) -- (10.23, 1.95);  %bf
\node at (10.56, 1.7) {\small 4};
\draw[->] (14.9, 2.1) -- (14.1, 3.9);  %cd
\node[right] at (14.5, 3) {6}; 
\node[right] at (15, 2){$V_3$};
\draw[->] (14.9, 2.05) -- (11.24, 3.86); %ce
\node at (14.3, 2.35) {\small 2};
\draw[<-] (14.8, 2) -- (10.3, 2); %cf
\node at (13.3, 2) {\small 5};
\draw[<-] (13.8, 4) -- (11.2, 4); %de
\node[above] at (12.5, 4) {4};
\node[above] at (14, 4) {$V_4$}; 
\draw[<-] (13.75, 3.85) -- (10.2, 2.1);  %df
\node at (13.05, 3.45)   {\small 11};
\node[above] at (11, 4) {$V_5$}; 
\draw[<-]  (10.83, 3.9) -- (10.11, 2.1); %ef
\node[left] at (10.5, 3) {7};

\node at (7.5, 2) {=};

\node[below] at (12, -1) {{\bf b.} $\tilde{\mathcal G}_{A, cpi}^6$};

%%%%%%%%%%%%%

\draw[<-] (21.18, 0) -- (23.9, 0); %ab
\node[below] at (21, 0) {$V_1$};
\node[below] at (22.5, 0) {1};
\draw[<-] (21.2, .07) -- (24.9, 2);   %ac
\node at (21.9, .5) {\small 2};
\draw[->] (21.15, .15) -- (23.9, 3.9); %ad
\node at (22.1, 1.5)  {\small 3};
\draw[->] (21, .2) -- (21, 3.9);  %a, e)
\node at  (21, 3) {\small 2}; 
\node[left] at (20, 2) {$V_6$}; 
\draw[dashed, <-] (20.9, .1) -- (20.1, 2);  %af
\node[left] at (20.5, 1) {2};
\draw[<-] (24.1, .1) -- (24.9, 1.9); %bc
\node[right] at  (24.5, 1) {6}; 
\node[below] at (24, 0) {$V_2$};
\draw[->] (24, .1) -- (24, 3.85);  %bd
\node at (24, 1) {\small 10};
\draw[<-]  (23.9, .1) -- (21.1, 3.9); %be
\node at (22.1, 2.5) {\small 4};
\draw[<-] (23.82, .07) -- (20.23, 1.95);  %bf
\node at (20.56, 1.7) {\small 1};
\draw[dashed, <-] (25, 2.1) -- (24.1, 3.9);  %cd
\node[right] at (24.5, 3) {6}; 
\node[right] at (25, 2){$V_3$};
\draw[<-] (24.9, 2.05) -- (21.24, 3.86); %ce
\node at (24.3, 2.35) {\small 7};
\draw[dashed, ->] (24.8, 2) -- (20.3, 2); %cf
\node at (23.3, 2) {\small 5};
\draw[dashed, ->]  (23.8, 4) -- (21.2, 4); %de
\node[above] at (22.5, 4) {4};
\node[above] at (24, 4) {$V_4$}; 
\draw[->] (23.75, 3.85) -- (20.2, 2.1);  %df
\node at (23.05, 3.4)   {\small 3};
\node[above] at (21, 4) {$V_5$}; 
\draw[<-]  (20.83, 3.9) -- (20.11, 2.1); %ef
\node[left] at (20.5, 3) {5};

\node at (17.5, 2) {+};

\node[below] at (22, -1) {{\bf c.} $\tilde{\mathcal G}_{A, cyclic}^6$};

 \node[below] at (12.5, -2) {{\bf Figure 5.} An incomplete $\mathcal G_A^6$};

\end{tikzpicture}

\subsection{Incomplete graphs}\label{sect: incomplete}

Incomplete graphs are similarly  reduced.  The approach is illustrated      with  
    $\mathcal G_A^6$ (Fig.~5a) where arcs $\widehat{V_1V_6}$, $\widehat{V_3V_4}$, $\widehat{V_3V_6}$, and $\widehat{V_4V_5}$ are  excluded.  To complete $\mathcal G^6_A$, add in  
     the missing arcs  with arbitrarily selected  lengths.    (For   Fig.~5a,  
  arcs of zero length were added.)   
   Denote the completed graph by  
     $\tilde{\mathcal G}_A^6$,   
   and compute $\tilde{\mathcal G}^6_{A, cpi}$ (Fig.~5b) and $\tilde{\mathcal G}^6_{A, cyclic}$ (Fig.~5c).  Closed path lengths in $\mathcal G^6_A$ agree with their lengths in  $\tilde{\mathcal G}^6_{A, cyclic}$.   (The four dashed arrows in Fig.~5c denote the  forbidden arcs.)  As the path $V_2 \stackrel{10}{\longrightarrow} V_4 \stackrel{3}{\longrightarrow}  V_6. \stackrel{5}{\longrightarrow} V_5 \stackrel{7}{\longrightarrow} V_3 \stackrel{2}{\longrightarrow} V_1 \stackrel{-1}{\longrightarrow} V_2$ uses most of the  largest allowed  $\tilde{\mathcal G}^6_{A, cyclic}$ leg lengths,  it is easy to show  that  its length of 26 is the longest $\tilde{\mathcal G}^6_{A, cyclic}$ Hamiltonian path. Thus, this   
     is the longest $\mathcal G^6_A$ Hamiltonian path with  
      the same length.

\begin{Thm}\label{thm: incomplete}  For an  incomplete  $\mathcal G_A^n$,  replace all non-admissible arcs  
 with arcs of arbitrary lengths to define $\tilde{\mathcal G}_A^n$;  compute $\tilde{\mathcal G}_{A, cyclic}^n$. The length of  a closed path   in  $\mathcal G_A^n$  equals  its length  in  $\tilde{\mathcal G}_{A, cyclic}^n$.
 
 The length of an admissible path starting at $V_j$ and ending at $V_k$  in the incomplete $\mathcal G^n_A$  is its length in $\tilde{\mathcal G}^n_{A, cyclic}$ plus the length of  the $V_j\to V_k$ arc in $\tilde{\mathcal G}^n_{A, cpi}$.   
\end{Thm}

 The concluding  Thm.~\ref{thm: incomplete} assertion allows  
 the $V_j\to V_k$ arc in $\tilde{\mathcal G}^n_{A, cpi}$ to   be  a forbidden  $\mathcal G^n_A$ choice.  As a Fig.~5  example, the $\mathcal G^6_A$ path $V_6  \stackrel{12}{\longrightarrow}  V_5  \stackrel{5}{\longrightarrow} V_3$ has length 17.  In $\tilde{\mathcal G}^6_{A, cyclic}$ this path  $V_6  \stackrel{5}{\longrightarrow}  V_5  \stackrel{7}{\longrightarrow} V_3$ has length 12, which is   added to the 5 length of the banned  $V_6 \stackrel{5}{\longrightarrow}  V_3$ arc in $\tilde{\mathcal G}^6_{A, cpi}$.

 \noindent{\em Proof:}  A  closed path's length in $\mathcal G_A^n$ is the same in $\tilde{\mathcal G}_A^n$ and (by Thm.~\ref{thm: gen decomp})  in  $\tilde{\mathcal G}_{A, cyclic}^n$. 
 
 A connected path's length in $\mathcal G^n_A$ is the same  in $\tilde{\mathcal G}^n_A$, which equals the sum of its lengths in $\tilde{\mathcal G}^n_{A, cpi}$ and $\mathcal G^n_{A, cyclic}$.  Its $\tilde{\mathcal G}^n_{A, cpi}$  length is that of its $V_j\to V_k$ arc,  which completes the proof.  $\square$  \smallskip

\subsection{Lower degrees of freedom}\label{sect: A lower}

 The    $\mathbb C^n_{A}$   graphs    
  have all (and only)  path information, which  
 simplifies deriving   closed path properties, computing lengths, and designing algorithms by using  the (somewhat predictive)  algebra of three-cycles.   
 General   path properties follow   from the  
 $\mathbb C^n_{A}$ basis, which characterizes all $\mathcal G^n_{A, cyclic}$ choices. 
 By seeking  
  general closed path properties, rather than just Hamiltonian circuits,  $\mathbb C^n_A$ is a best possible component.  This is because the simplest closed path is   a triplet, so  
  that  three-cycle must be in $\mathbb C^n_A$.

\section{Symmetric cost  settings}

Deriving the  structure of $\mathbb G^n_S$---the space of  n-vertex complete symmetric weighted (no loops)  graphs---mimics
what was done    for $\mathbb G^n_A$.  As with Eq.~\ref{eq: first decomp},  the goal is to decompose  a $\mathcal G^n_S\in \mathbb G^n_S$ as \begin{equation}\label{eq: second decomp}  \mathcal G^n_S = \mathcal G^n_{S, cpi} + \mathcal G^n_{S, cyclic}, \end{equation} where the definition of $ \mathcal G^n_{S, cpi}$  will captures terms that, at least initially, can be ignored   when computing lengths of $\mathcal G^n_S$ closed paths.   
  Thus,  all central  $\mathcal G^n_S$ path properties  are based on the    $ \mathcal G^n_{S, cyclic}$ structure.   
The approach   is to find   
 the manifold  
 of  $\mathbb G^n_S$ graphs where closed paths have a fixed length; this manifold defines the $\mathcal G^n_{S, cpi}$ components.  The    normal bundle of this manifold  measures   deviations from ``sameness"  to capture what is needed to find  closed  path properties of a  $\mathcal G_S^n \in \mathbb G^n_S$.  

     \begin{tikzpicture}[xscale=0.56, yscale=0.48]
\draw (-7.9, 0) -- (-5.1, 0);  %ab
\draw (-5, .1) --(-5, 2.9);  % bc
\draw (-5.1, 3) -- (-7.9, 3); %cd
\draw (-8, .1) -- (-8, 2.9);   %ad
\draw  (-5.1, .1) -- (-7.9, 2.9);  %bd
\draw (-7.9, .1) -- (-5.1, 2.9);  % ac
 
\node[below] at (-6.5, 0) {$y_1$};
\node[right] at (-5, 1.5) {$x_2$};
\node[above] at (-6.5, 3) {$y_2$};
\node[left] at (-8, 1.5) {$x_1$};
\node[below] at (-7, 2.5) {$z_2$};
\node[below] at (-5.75, 2.5) {$z_1$};
 
\node[below] at (-8.15, 0.15) {$V_1$};
\node[below] at (-4.75, 0.15) {$V_2$};
\node[above] at (-8.15, 2.9) {$V_4$};
\node[above] at (-4.75, 2.9) {$V_3$};

\node[below] at (-6.5, -.9) {{\bf a.} ${\mathcal G}_S^4$};

%%%%%%%%%%%%%%%%%

\node [below] at (1.1, 3.6) {$V_1  \stackrel{y_1}{\longrightarrow} V_2  \stackrel{x_2}{\longrightarrow} V_3  \stackrel{y_2}{\longrightarrow} V_4  \stackrel{x_1}{\longrightarrow} V_1$};

\node [below] at (1.1, 2.3) {$V_1  \stackrel{z_1}{\longrightarrow} V_3  \stackrel{y_2}{\longrightarrow} V_4  \stackrel{z_2}{\longrightarrow} V_2  \stackrel{y_1}{\longrightarrow} V_1$};

\node [below] at (1.1, 1) {$V_1  \stackrel{x_1}{\longrightarrow} V_4  \stackrel{z_1}{\longrightarrow} V_2  \stackrel{x_2}{\longrightarrow} V_3  \stackrel{z_2}{\longrightarrow} V_1$};

\node[below] at (1.5, -.9) {{\bf b.} Equal path lengths};

%%%%%%%%%%

\draw (8.1, 0) -- (10.9, 0);
\draw (11, .1) --(11, 2.9);
\draw (10.9, 3) -- (8.1, 3);
\draw (8, .1) -- (8, 2.9);
\draw (10.9, .1) -- (8.1, 2.9);
\draw (8.1, .1) -- (10.9, 2.9);
\node[below] at (9.5, 0) {$\scriptstyle{\omega_1+\omega_2} $};
\node[right] at (11, 1.5) {$\scriptstyle{\omega_2+\omega_3}$};
\node[above] at (9.5, 3) {$\scriptstyle{\omega_3+\omega_4}$};
\node[left] at (8, 1.5) {$ \scriptstyle{\omega_1 +\omega_4}$};
\node at (9, 2.3) {$\scriptstyle{\omega_2+\omega_4}$};
\node at (9, 1.0) {$\scriptstyle\scriptstyle{\omega_1+\omega_3}$};
\node[below] at (7.55, 0.15) {$ \omega_1$};
\node[below] at (11.35, 0.15) {$ \omega_2$};
\node[above] at (7.55, 2.9) {$ \omega_4$};
\node[above] at (11.35, 2.9) {$ \omega_3$};
\node[below] at (9.5, -.9) {{\bf c.} $\mathcal G^4_{S, cpi}$ representation};

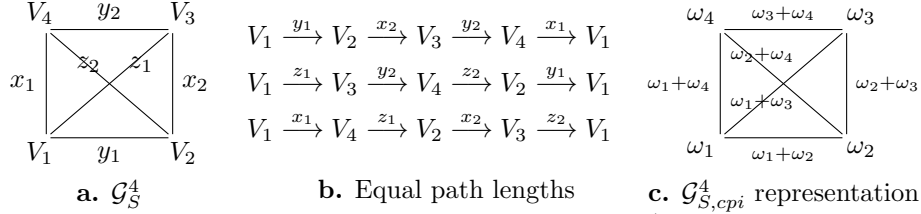
\node[below] at (.5, -1.85) {{\bf Figure 6.} Closed path independence for $\mathbb G^4_S$};

\end{tikzpicture}

The  $\mathbb G^n_S$  cpi definition and structure  differs   from that of  $\mathbb G^n_A$; e.g.,  
 rectangles replace triangles.  For instance, $\mathcal G_S^4 \in \mathbb G^4_S$ (Fig.~6a) is cpi iff   
the three  Fig.~6b  routes have equal length.  With cancelations, this requires the three sums of the vertical, the  horizontal,   and  the  diagonal  lengths to agree, or   \begin{equation}\label{eq: s 4-path} x_1+x_2=y_1+y_2 = z_1+z_2.\end{equation}  
  Hamiltonian paths for  $\mathcal G_S^4\in 
 \mathbb G^4_S$   combine two of the pairs of diagonals,   vertical edges,  or horizontal edges, so  the two smallest Eq.~\ref{eq: s 4-path} sums define the shortest  path with  length  equal to  this sum.\footnote{Should the vertices define a triangle with one in the interior, the three pairs are defined by the triangle's three vertices.  A pair is the arc from a vertex to the interior point and the  triangle's leg that is opposite the vertex.} For instance, if $x_1+x_2=10, \,    y_1+y_2=20,$ and $z_1+z_2 = 30$, then the shortest Hamiltonian circuit traverses the perimeter and it has length 30..

All cpi graphs in  $\mathbb G^4_S$   
  satisfy    two independent equations (Eq.~\ref{eq: s 4-path}) in  six variables.  One solution has zero leg lengths,     
 so  all solutions (i.e., all  cpi graphs $\mathcal G^4_{S, cpi}$)  are characterized by  Eq.~\ref{eq: s 4-path}'s four-dimensional kernel.   
One choice  uses 
 weights $\{\omega_j\}_{j=1}^4$ where $\omega_j$  is assigned to vertex $V_j$, $j=1, \dots, 4$,  to define  
 the $\widehat{V_jV_k}$  length of ${\omega_j+\omega_k}$  (Fig.~6c).  It follows from Fig.~6c that this choice  satisfies  Eq.~\ref{eq: s 4-path}.   For $\mathbb G^4_S$,   these are the   
`closed path independent' graphs.  The common path  length   depends on how often each vertex is visited; e.g., a    $\mathcal G^4_{S, cpi}$  closed path that visits each of the three vertices $\{V_i\}_{i=1}^3$ twice has length $4\sum_{j=1}^3 \omega_j$;    Hamiltonian paths in $\mathcal G^4_{S, cpi}$  have  length  $T(\mathcal G_{S, cpi}^4) = 2\sum_{j=1}^4 \omega_j$.

The above  discussion centered on  Fig.~6  extends to  $n\ge4.$  
 
\begin{Def}\label{def: S cpi}  A graph $\mathcal G^n_S \in \mathbb G^n_S$ is {\em `closed path independent'}  iff for any set of vertices, all closed paths that pass through each of these vertices once have the same length.\end{Def}  
  
\begin{Thm}\label{prop: omega} For $n\ge 4$, a cpi    graph  $\mathcal G^n_{S, cpi} \in \mathbb G^n_S$  assigns
  weight $\omega_j$  to vertex $V_j$, $j=1, \dots, n;$  the 
 $\widehat{V_jV_k}$  length is  ${\omega_j+\omega_k}$, $j\ne k$. A closed path passing once through the vertices $\{V_j\}_{j\in \mathcal D}$ has length  $2\sum_{j\in \mathcal D} \omega_j$.  A Hamiltonian path length in $\mathcal G^n_{S, cpi}$  is $2\sum_{j=1}^n \omega_j$.   
 \end{Thm}

 Determining   
  the structure of    $\mathcal G^n_{S, cpi}$  
  requires identifying  
  $\mathbb G^n_S$ with $\mathbb R^{n\choose2}_S$.  Here, $\mathbb R^{n\choose2}_S$ differs from   $\mathbb R^{n\choose2}_A$ (Eq.~\ref{def: dn})   
  because  in  $\mathbb G^n_A$  (Sect.~\ref{sect: compute}),   
 $d_{i, j}=-d_{j, i}$, but  in   
  $\mathbb G^n_S$,  $d_{i, j}=d_{j, i}$.  Thus,  \begin{equation}\label{def: dn S}\mathbf d_S^n=(d_{1, 2}, d_{1, 3}, \dots, d_{1, n}; d_{2, 3}, \dots, d_{2, n}; d_{3, 4}, \dots ; d_{n-1, n}) \in \mathbb R^{n\choose2}_S, \textrm{ where } d_{i, j} = d_{j, i};\end{equation}

\begin{Thm}\label{thm: symmetric cpi}  The space of $\mathcal G^n_{S, cpi}$ graphs, denoted by $\mathbb {CPI}_S^n$,  is a $n$-dimensional linear subspace of $\mathbb G^n_S$, or, equivalently, of $\mathbb R^{n\choose2}_S.$  
Let $\mathbf B^n_j \in \mathbb R_S^{n\choose2}$  be where $d_{j, k} = 1$ for all $k\ne j, k=1, \dots, n$; all other $d_{u, v}=0$.  A basis for $\mathbb {CPI}_S^n$ is  $\{\mathbf B^n_j\}_{j=1}^n$.  
\end{Thm}

As $\mathcal G^n_{S, cpi}$  will identify   
 components of $\mathcal G_S^n$  entries that cloud analyzing       closed paths, the  
 emphasis shifts  to    $\mathcal G^n_{S, cyclic} = \mathcal G_S^n-\mathcal G^n_{S, cpi}$.  The dimensions of    $\mathbb G_S^n$ and   $\mathbb{CPI}_S^n$   are    $n\choose2$ and   
   $n$,  so $\mathbb{CPI}_S^n$'s normal subspace, $\mathbb C^n_S$,  has dimension $\frac{n(n-3)}2$. (The dimension of $\mathbb C^n_A$ is ${{n-1}\choose2} = \frac{n(n-3)}2 + 1$.) Its four-cycle  structure (Thm~\ref{thm: s cyclic basis})  
   identifies    $\mathcal G^n_S$'s inherent symmetry.   

\begin{Thm}\label{thm: s cyclic basis}  Let vector $\mathbf b^n_{i, j, k, s}$ be  where $d_{i, j} = 1, d_{j, k}=-1, d_{k, s}=1, d_{s, i}=-1$; all other $d_{u, v}=0$.   
  The space $\mathbb C^n_S$  spanned by all $\{\mathbf b^n_{i, j, k, s}\}$, with dimension $ \frac{n(n-3)}2$, is orthogonal to $\mathbb{CPI}_S^n$.  With  
\begin{equation}\label{eq: CS basis}\mathcal A^n_{1,2} =  \{\mathbf b^n_{1, 2, j, k}\}_{2<j<k\le n}, \,  \mathcal B^n_{1, 2} =  \{\mathbf b^n_{1, 2, j, 3}\}_{j=4}^n, \textrm{ a basis is } \mathcal A^n_{1, 2} \cup \mathcal B^n_{1, 2}.\end{equation}  \end{Thm}

 To explain these   $\mathbf b^n_{i, j, k, s}$   vectors,   if the arc lengths of the route around the perimeter of Fig.~6a  satisfy the  Eq.~\ref{eq: s 4-path} cpi requirement, then  
  $x_1-y_1+x_2-y_2=0$.  This equation has the normal vector $(1, -1, 1, -1)$, which in path form, is  $\mathbf b^4_{4, 1, 2, 3}$.   In general the 
 $\mathbf b^n_{i, j, k, s}$ path   is the four-cycle $V_i \stackrel{1}{\longrightarrow} V_j \stackrel{-1}{\longrightarrow} V_k \stackrel{1}{\longrightarrow} V_s \stackrel{-1}{\longrightarrow} V_i$;  each    vertex   has one leg of length $1$  and one of length $-1$.    With  Fig.~6a,  the   $\mathbf b^4_{1, 2, 3, 4}$ and $\mathbf b^4_{1, 4, 2, 3}$ multiples are, respectively,  $\frac14\{(y_1+y_2) -(x_1+x_2)\}$  and  $\frac14\{(x_1+x_2) - (z_1+z_2)\}$.     Namely,   $\mathbf b^n_{i, j, k, s}$ measures how the $\mathcal G^n_S$  data deviates   from   cpi sameness while identifying which data edges of a four-tuple   are    ridges or valleys.      \medskip

\noindent{\em Proof.} That  $\mathbf b^n_{i, j, k, s}$  is orthogonal to each $\mathbf B^n_t$ is immediate.   If $t\ne i, j, k, s$,  the scalar product is zero.  If $t$ is one of these indices, say $t=j$, then one component of  $\mathbf b^n_{i, j, k, s}$ with vertex $V_j$ is positive and the other is negative, so the scalar product with $\mathbf B^n_t$ is zero. 

Establishing  the linear independence of  Eq.~\ref{eq: CS basis} follows a switching  pattern.  Iteratively, it will be shown that all  coefficients of $\sum_{2<j<k\le n} \alpha_{j, k} \mathbf b^n_{1, 2, j, k} = \mathbf 0$ must equal zero.  For each of the ${n-3}\choose2$  top vectors in $\mathcal A^n_{1, 2}$ (i.e., $j, k\ge 4$), only vector $\mathbf b^n_{1, 2, j, k}$ has a non-zero $d_{j, k}$, so $\alpha_{j, k}=0.$  For all remaining vectors,   either $j$ or $k$ equals 3.  Of these,  only  
the top  $\mathcal B^n_{1, 2}$ vector of  $\mathbf b^n_{1, 2, n, 3}$  has a non-zero $d_{2, n}$, so $\alpha_{n, 3}=0$. The top remaining $\mathcal A^n_{1, 2}$ vector  is $\mathbf b^n_{1, 2, 3, n}$, where,  with the removal of  $\mathbf b^n_{1, 2, n, 3}$,  only $\mathbf b^n_{1, 2, 3, n}$ has non-zero $d_{3, n}$, so $\alpha_{3, n}=0$.  
  The obvious induction argument of switching between remaining  $\mathcal A^n_{1, 2}$ and $\mathcal B^n_{1, 2}$ vectors continues.  
That is,   if $s$ is the upper bound of the remaining $j, k$ values, then 
  only $\mathbf b^n_{1, 2, s, 3} \in \mathcal B^n_{1, 2}$ has a non-zero $d_{2, s}$ term, so $\alpha_{s, 3}=0$.  The top remaining  $\mathcal A^n_{1, 2}$ vector is $\mathbf b^n_{1, 2, 3, s}$; as   $\mathbf b^n_{1, 2, s, 3}$ was removed,  
   only $\mathbf b^n_{1, 2, 3, s}$ of the remaining vectors has a non-zero  $d_{3, s}$ term, so $\alpha_{3, s}=0$  and $s-1$   is the largest  remaining $j, k$ value.  This completes the proof. $\square$\medskip

\subsection{Decomposing  $\mathbb G^n_S$}

Theorem \ref{thm:  GS structure} summarizes  the above; it is the    $\mathbb G^n_S$ version of  Thm.~2.   As with Sect.~2, the decomposition involves $O(n^2)$ computations. 
\begin{Thm} \label{thm:  GS structure} For $n\ge 4$,    $\mathbb G^n_S$ has an $n$-dimensional linear subspace $\mathbb {CPI}_S^n$ and an orthogonal $\frac{n(n-3)}2$ dimensional linear subspace $\mathbb C^n_S$.   
A $\mathcal G_S^n\in \mathbb G^n_S$ has a unique decomposition $ 
  \mathcal G_S^n = \mathcal G^n_{S, cpi} + \mathcal G^n_{S, cyclic}$  
  where $\mathcal G^n_{S, cpi}\in \mathbb{CPI}^n_S$ and $\mathcal G^n_{S, cyclic} \in \mathbb C^n_S$  are, respectively, the orthogonal projection of $\mathcal G^n_S$ to  $\mathbb{CPI}^n_S$ and to $\mathbb C^n_S$. \end{Thm}  

\begin{tikzpicture}[xscale=0.45, yscale=0.45]

\draw (0, 0) -- (2.9, 0);
\draw  (0, 0) -- (0, 2.9);
\draw (3, 0) -- (3, 2.9); 
\draw (0.1, 3) -- (2.9, 3);
\draw (0.1, 2.9) -- (2.9, 0.1);
\draw (0.1, 0.1) -- (2.9, 2.8);
\node[below] at (0, 0) {$V_1$};
\node[below] at (3, 0) {$V_2$};
\node[above] at (3, 3) {$V_3$};
\node[above] at (0, 3) {$V_4$}; 
\node[below] at (1.5, 0) {17};
\node[left] at (0, 1.5) {13};
\node[right] at (3, 1.5) {17};
\node[above] at (1.5, 3) {17};
\node at (0.9, 1) {$14$};
\node at (0.9, 2) {$12$};

\node[below] at (1.5, -1) {{\bf a.} $\mathcal G^4_S$};

\node at (4.5, 1.5) {=};

%%%%%%%%%%%%%

\draw (6, 0) -- (8.9, 0);  %ab
\draw  (6, 0) -- (6, 2.9); %ad
\draw (9, 0) -- (9, 2.9);  %bc
\draw (6.1, 3) -- (8.9, 3); %
\draw (6.1, 2.9) -- (8.9, 0.1);
\draw (6.1, 0.1) -- (8.9, 2.8);
\node[below] at (6, 0) {$V_1$};
\node[below] at (9, 0) {$V_2$};
\node[above] at (9, 3) {$V_3$};
\node[above] at (6, 3) {$V_4$}; 
\node[below] at (7.5, 0) {15};
\node[left] at (6, 1.5) {13};
\node[right] at (9, 1.5) {17};
\node[above] at (7.5, 3) {15};

\node at (6.9, 1) {$16$};
\node at (6.9, 2) {$14$};

\node[below] at (7.5, -1) {{\bf b.} $\mathcal G^4_{S, cpi}$};

\node at (10.5, 1.5) {+};

%%%%%%%%%%%%%%%%%

\draw (12, 0) -- (14.9, 0);  %ab
\node[below] at (13.5, 0) {2};
\draw  (12, 0) -- (12, 2.9); %ad
%\node[left] 
\draw (15, 0) -- (15, 2.9);  %bc
\draw (12.1, 3) -- (14.9, 3); %cd

\draw (12.1, 2.9) -- (14.9, 0.1);
\draw (12.1, 0.1) -- (14.9, 2.8);
\node[below] at (12, 0) {$V_1$};
\node[below] at (15, 0) {$V_2$};
\node[above] at (15, 3) {$V_3$};
\node[above] at (12, 3) {$V_4$}; 
 
\node[left] at (12, 1.5) {0}; %ad
\node[right] at (15, 1.5) {0}; %bc
\node[above] at (13.5, 3) {2};  %cd
\node at (12.9, 1) {$-2$};
\node at (12.9, 2) {$-2$}; 

\node[below] at (13.5, -1) {{\bf c.} $\mathcal G^4_{S, cyclic}; 2\mathbf b_{1, 2, 4, 3}$};

\node at (10.5, 1.5) {+};

\node at (7.5, -2.8) {{\bf Figure 7.} Interpreting $\mathcal G^n_{S, cyclic}$};

\end{tikzpicture}

 Before      computing  $\mathcal G^n_{S, cpi}$ and $\mathcal G^n_{S, cyclic}$,  Fig.~7 is used to explain  
  their roles  and to relate $\mathcal G^n_{S, cpi}$ with $\mathcal G^n_{A, cpi}.$   According to Thm.~\ref{thm:  GS structure},   $\mathcal G^n_{S, cpi}$ is the $\mathbb{CPI}_S^n$ graph that most closely resembles $\mathcal G^n_S$; this  
 similarity   is apparent when comparing  Figs.~7a, b. A defining feature of $\mathcal G^n_{A, cpi} \in \mathbb{ST}_A^n$ is that for any set of vertices, the length of all closed paths meeting each vertex once is zero.   
 Similarly, according to Def.~\ref{def: S cpi}, for $\mathcal G^n_{S, cpi} \in \mathbb {CPI}_S^n$ and for any selected set of vertices,  the length of all closed paths that meet each  vertex once is the same, but  not necessarily zero. The role of the decomposition in $\mathbb G^n_A$ and in $\mathbb G^n_S$ is to remove these common path length values.    
  Thus the $\mathcal G^n_{A, cyclic}$ and $\mathcal G^n_{S, cyclic}$ graphs characterize how the data from the original graph ($\mathcal G^n_A$ or $\mathcal G^n_S$) differs from the cpi sameness to  provide   valued path length information.

In  
  Fig.~7b,    
  the sums of its vertical edges,  horizontal edges, and  diagonals all equal  30. Thus, all  Fig.~7b  Hamiltonian paths have the  length 60.  
  Three of the six
    Fig.~7a  Hamiltonian circuits  are $V_1 \stackrel{17}{\longrightarrow} V_2 \stackrel{17}{\longrightarrow} V_3 \stackrel{17}{\longrightarrow} V_4 \stackrel{13}{\longrightarrow} V_1$ with length 64, $V_1 \stackrel{17}{\longrightarrow} V_2 \stackrel{12}{\longrightarrow} V_4  \stackrel{17}{\longrightarrow} V_3 \stackrel{14}{\longrightarrow} V_1$ with length 60, and  $V_1 \stackrel{14}{\longrightarrow} V_3 \stackrel{17}{\longrightarrow} V_2  \stackrel{12}{\longrightarrow} V_4 \stackrel{13}{\longrightarrow} V_1$ with length 56; the other three are reversals.  The average  
    length  of these  paths
   is 60, which agrees with its Fig.~7b value. This comparison  accurately suggests that for any set of vertices used to define  
     closed paths, what happens in  $\mathcal G^n_{S, cpi}$  is the average of what happens in $\mathcal G^n_S.$   For instance, the length of a  Hamiltonian path for Fig.~7a  is  the sum of its lengths in Fig.7b and Fig.~7c.  A ``subtraction" argument, similar to that used with Fig.~2, is that the  portion of a path entry contributing to the Fig.~7b  average Hamiltonian length is subtracted from the actual leg value.  What  remains  determines how the path length differs from the average, so it is used in the computation; the average length of a Hamiltonian path is replaced at the end. 
 Thus, as developed below,  
   path lengths in $\mathcal G^n_{S, cyclic}$ (e.g.,  Fig.~7c) measure differences  
    from the average.   As  $V_1 \stackrel{-2}{\longrightarrow} V_3 \stackrel{0}{\longrightarrow} V_2  \stackrel{-2}{\longrightarrow} V_4 \stackrel{0}{\longrightarrow} V_1$  in Fig.~7c has the shortest length of $-4$,    
   this  defines  the shortest  Fig.~7a path that has length $-4$ from the average of 60, or 56.

\subsection{Computing  $\mathcal G^n_{S, cpi}$ and $\mathcal G^n_{S, cyclic}$ }
 
Computing the $\mathcal G^n_{S, cpi}$ and $\mathcal G^n_{S, cyclic}$ for a  
 $\mathcal G^n_S$ follows the lead of Sect.~\ref{sect: compute}.  This is because $\mathcal G^n_{S, cpi}$  is the orthogonal projection of $\mathcal G^n_S$ to $\mathbb{CPI}^n_S$, and a basis for $\mathbb{CPI}^n_S$ is known (Thm.~\ref{thm: symmetric cpi}).
 Entries for $\mathcal G^n_{S, cpi}$ and $\mathcal G^n_{S, cyclic}$ are based on the following.    
   
\begin{Def}\label{def: S for SS} For     $V_j$  
 of $\mathcal G_S^n \in \mathbb G^n_S$,     let   $\mathcal S_S(V_j)$ be  the sum of the arc lengths attached to vertex $V_j$, $j=1, \dots, n.$   Let $T(\mathcal G_{S}^n) = \frac1{n-1}  \sum_{j=1}^n \mathcal S_S(V_j)$. \end{Def}

Because $\frac1{n-1}S(V_j)$  is the average length of an arc with $V_j$ as a vertex, it follows that $T(\mathcal G_{S}^n)$ is the {\em    average $\mathcal G^n_S$ Hamiltonian path length.}    
  As $\mathcal G^n_{S, cyclic}$  consists of   $\mathbf b^n_{i, j, k, s}$ cycles,     each   $\mathbf b^n_{i, j, k, s}$ arc entering   a vertex has  
 a leaving arc with the same weight but opposite sign, so $S_S(V_j)=0$.   
 This equation requires 
 the $S_S(V_j)$ values for $\mathcal G^n_{S}$ and $\mathcal G^n_{S, cpi}$ to agree.  Because  $T(\mathcal G^n_S)$ sums these values,   the average Hamiltonian path lengths in  $\mathcal G^n_{S}$ and in $\mathcal G^n_{S, cpi}$ agree (as suggested with Fig.~7), or   
 \begin{equation}\label{eq: t equals t} T(\mathcal G^n_{S}) = T(\mathcal G^n_{S, cpi}) = 2\sum_{j=1}^n \omega_j.\end{equation}

Agreement between  
   $S_S(V_j)$ values in  $\mathcal G_S$ and  $\mathcal G^n_{S, cpi}$  provides   equations  for the unknowns $\{\omega_j\}_{j=1}^n$.   
    Illustrating with  Fig.~8a,  as  $S_S(V_1) =  81$ for $\mathcal G^5_S$,  
 the same value holds for $\mathcal G^5_{S, cpi}$, which means that
$\sum_{j=2}^5( \omega_1+\omega_j)=81$.  
In general, the unknown  $\{\omega_j\}_{j=1}^n$  satisfy  
 \begin{equation}\label{eq: leading to omega} S_S(V_j) = \sum_{k\ne j} (\omega_j+\omega_k) = (n-1)\omega_j + \sum_{k\ne j} \omega_k =  (n-2)\omega_j + \sum_{k=1}^n \omega_k = (n-2)\omega_j + \frac12T(\mathcal G^n_{S, cpi}),.\end{equation} 
 Using $T(\mathcal G^n_S) = \mathcal T(\mathcal G^n_{S, cpi})$  (Eq.~\ref{eq: t equals t}),  
the values of the $\mathcal G^n_{S, cpi}$ weights are   
   \begin{equation}\label{def: w eq}  \omega_j =  
    \frac{1}{n-2}[S_S(V_j) - \frac12T(\mathcal G_{S}^n)], \,  j=1,  2, \dots, n. \end{equation}

These $\omega_j$ weights, which define  
 $\mathcal G^n_{S, cpi}$ and  $\mathcal G^n_{S, cyclic}$,   lead to a result concerning path lengths.

\begin{Thm}\label{thm: 10 done better}   
For $\mathcal G^n_S\in \mathbb G^n_S$,  
Eq.~ \ref{def: w eq} defines the weights of its
  $\mathcal G^n_{S, cpi}$ component.  Let $\mathcal G^n_{S, cyclic} = \mathcal G^n_S-\mathcal G^n_{S, cpi}.$   
The  $\mathcal G_S^n$  length of a Hamiltonian circuit  equals $T(\mathcal G_S^n)$ plus its  $\mathcal G^n_{S, cyclic}$ path length. \end{Thm}  

\noindent{\em Proof:} A  $\mathcal G^n_S$ path length   is the sum of its   $\mathcal G^n_{S, cpi}$ and   $\mathcal G^n_{S, cpi}$ lengths.   
All  $\mathcal G^n_{S, cpi}$ Hamiltonian paths  have  length $T(\mathcal G^n_{S, cpi}) = T(\mathcal G_S^n)$, so Thm.~\ref{thm: 10 done better} follows. $\square$     

According to Thm.~\ref{thm: 10 done better}, all essential closed path properties of $\mathcal G^n_S$  are based on  the structure of  $\mathcal G^n_{S, cyclic}$ and its four-cycle symmetries.  Thus, general properties characterizing  $\mathcal G^n_{S, cyclic}$ are useful.  
\begin{Cor} \label{cr: cyclic zero} If $\mathcal G^n_S$ has the property that $S_S(V_j)=0$, $j=1, \dots, n$, then $\mathcal G^n_S\in \mathbb C^n_S$.\end{Cor}

{\em Proof:}  
This condition requires $T(\mathcal G^n_S)=0$ (Def.~\ref{def: S for SS}) and $\omega_j=0$, $j=1, \dots, n$ (Eq.~\ref{def: w eq}).  As $\mathcal G^n_{S, cpi}= 0$, it follows that $\mathcal G^n_S = \mathcal G^n_{S, cyclic}$. $\square$\smallskip 

 \begin{tikzpicture}[xscale=0.5, yscale=0.5]

 \node[below] at (1, 0) {$V_1$};
 
\node[below] at (4, 0) {$V_2$};
 
\node[right] at (5, 2) {$V_3$};
 
\node[above] at (2.5, 3) {$V_4$};
 
\node[left] at (0, 2) {$V_5$};
\draw  (1.1,  0) -- (3.8, 0); %ab
\draw   (4.9, 1.9) -- (1.15, .15);  %ca
\draw    (1, 0) -- (2.4, 2.9); %ad
\draw   (.9, .05) -- (.1, 1.9); %ae
\draw    (0.15, 2) -- (3.8, .15);  %be
 \draw   (2.55, 3) -- (3.9, .15);  %db
 \draw    (4, 0) -- (4.9, 1.8); %bc
\draw    (4.9, 2.1) -- (2.75, 3); %cd
\draw   (4.85, 2) -- (0.45, 2); %ce
\draw   (0.2, 2.1) -- (2.3, 2.95); %ed

\node[below] at (2.4, 0) {$3$};  %ab
\node[right] at (4.5, 1) {$9$}; %bc
\node[left] at (0.5, 1) {$24$}; %ae
\node[above] at (3.75, 2.5) {$18$};  %cd
\node[above] at (1.25, 2.5) {$9$}; %ed

\node[right] at (2.52, 1) {\small 12}; %AC
 
\node[left] at (2.25, 1.7) {\small 42}; %ad
 
\node[right] at (1.5, 1) {\small 27};
 
\node[above] at (3.34, 1.3) {\small 21}; %bd
 
\node[left] at (2.9, 2.1) {\small 3}; %ec

\node[below] at (2.5, -1) {{\bf a.} Original; $\mathcal G_S^5$}; 
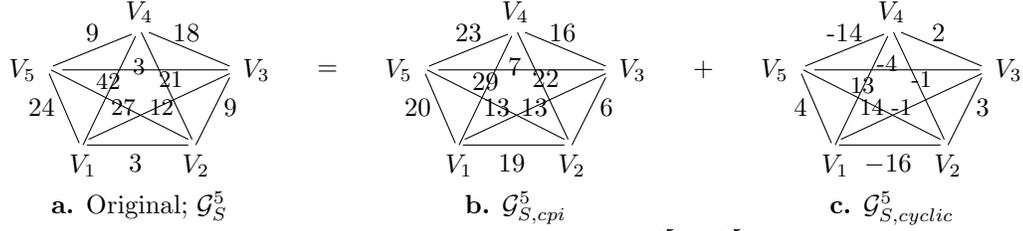
\node[below] at (12.5, -2) {{\bf Figure 8.} Decomposition of a $\mathcal G_S^5 \in \mathbb G^5_S$};
\node at (7.5, 2) {=};
\node at (17.5, 2) {+};

%%%%%%%%%%%%%%%

\node[below] at (11, 0) {$V_1$};
 
\node[below] at (14, 0) {$V_2$};

\node[right] at (15, 2) {$V_3$};
 
\node[above] at (12.5, 3) {$V_4$};
 
\node[left] at (10, 2) {$V_5$};
\draw (11.1,  0) -- (13.8, 0); %ab
\draw  (14.9, 1.9) -- (11.15, .15);  %ca
\draw   (11, 0) -- (12.4, 2.9); %ad
\draw  (10.9, .05) -- (10.1, 1.9); %ae
\draw   (10.15, 2) -- (13.8, .15);  %be
 \draw (12.55, 3) -- (13.9, .15);  %db
 \draw  (14, 0) -- (14.9, 1.8); %bc
\draw  (14.9, 2.1) -- (12.75, 3); %cd
\draw  (14.85, 2) -- (10.45, 2); %ce
\draw (10.2, 2.1) -- (12.3, 2.95); %ed

\node[below] at (12.4, 0) {$19$};  %ab
\node[right] at (14.5, 1) {$6$}; %bc
\node[left] at (10.5, 1) {$20$}; %ae
\node[above] at (13.75, 2.5) {$16$};  %cd
\node[above] at (11.25, 2.5) {$23$}; %ed

\node[right] at (12.4, 1) {$ 13$}; %AC
 
\node[left] at (12.3, 1.7) {$29$}; %ad
 
\node[right] at (11.4, 1) {$13$};
 
\node[above] at (13.3, 1.3) {$22$}; %bd
 
\node[left] at (12.9, 2.1) {$7$}; %ec

\node[below] at (12.5, -1) {{\bf b.} $\mathcal G^5_{S, cpi}$}; 

%%%%%%%%%%%%%%%%%%%%%
%ab=1, ac=-4, ad=-3, ae=-9, bc=-5, bd=-4, be=-10; cd=1, ce=-5, de=-6
\node[below] at (21, 0) {$V_1$};
 
\node[below] at (24, 0) {$V_2$};

\node[right] at (25, 2) {$V_3$};
 
\node[above] at (22.5, 3) {$V_4$};
 
\node[left] at (20, 2) {$V_5$};
\draw (21.1,  0) -- (23.8, 0); %ab
\draw  (24.9, 1.9) -- (21.15, .15);  %ca
\draw   (21, 0) -- (22.4, 2.9); %ad
\draw  (20.9, .05) -- (20.1, 1.9); %ae
\draw  (20.15, 2) -- (23.8, .15);  %be
 \draw  (22.55, 3) -- (23.9, .15);  %db
 \draw   (24, 0) -- (24.9, 1.8); %bc
\draw  (24.9, 2.1) -- (22.75, 3); %cd
\draw (24.85, 2) -- (20.45, 2); %ce
\draw  (20.2, 2.1) -- (22.3, 2.95); %ed

\node[below] at (22.4, 0) {$-{16}$};  %ab
\node[right] at (24.5, 1) {$3$}; %bc
\node[left] at (20.5, 1) {$4$}; %ae
\node[above] at (23.75, 2.5) {$2$};  %cd
\node[above] at (21.25, 2.5) {-14}; %ed

\node[right] at (22.25, 1) {\small-1}; %AC
 
\node[left] at (22.3, 1.6) {\small{13}}; %ad
 
\node[right] at (21.4, 1) {\small 14};
 
\node[above] at (23.3, 1.3) {\small -1}; %bd
 
\node[left] at (22.9, 2.18) {\small -4}; %ec

\node[below] at (22.5, -1) {{\bf c.} $\mathcal G^5_{S, cyclic}$};

\end{tikzpicture}

To illustrate   Thm.~\ref{thm: 10 done better},  
 the Fig.~8a computations  from $\mathcal G_S^5$    
are $S_S(V_1) = 81, S_S(V_2) = 60, S_S(V_3) = 42,  S_S(V_4) = 90, S_S(V_5) = 63$, so $T(\mathcal G_S^5) =  84.$   
This means that (Eq.~\ref{def: w eq})  $\omega_1 = \frac13[81-42]= 13, \, \omega_2 =  6, \,  \omega_3=  0, \,  \omega_4= 16, \, \omega_5 = 7,$ from which   $\mathcal G^5_{S, cpi}$ and $\mathcal G^5_{S, cyclic}$ of Figs.~8b, c  follow.

\begin{tikzpicture}[xscale=0.5, yscale=0.52]

\draw (1.18, 0) -- (3.9, 0); %ab
\node[below] at (1, 0) {$V_1$};
\node[below] at (2.5, 0) {14};
\draw (1.15, .06) -- (4.9, 2);   %ac
\node at (1.9, .5) {\small 11};
\draw (1.15, .15) -- (3.9, 3.9); %ad
\node at (2.1, 1.5)  {\small 15};
\draw (1, .2) -- (1, 3.9);  %a, e)
\node at  (1, 3) {\small 17}; 
\node[left] at (0, 2) {$V_6$}; 
\draw (.9, .1) -- (.1, 2);  %af
\node[left] at (.5, 1) {12};
\draw (4.1, .1) -- (4.9, 1.9); %bc
\node[right] at  (4.5, 1) {9}; 
\node[below] at (4, 0) {$V_2$};
\draw (4, .1) -- (4, 3.85);  %bd
\node at (4, 1) {\small 9};
\draw  (3.9, .1) -- (1.1, 3.9); %be
\node at (2.1, 2.5) {\small 15};
\draw (3.82, .07) -- (.23, 1.95);  %bf
\node at (.56, 1.7) {\small 10};
\draw (4.9, 2.1) -- (4.1, 3.9);  %cd
\node[right] at (4.5, 3) {13}; 
\node[right] at (5, 2){$V_3$};
\draw (4.9, 2.05) -- (1.24, 3.86); %ce
\node at (4.3, 2.35) {\small 19};
\draw (4.8, 2) -- (.2, 2); %cf
\node at (3.3, 2) {\small 13};
\draw (3.8, 4) -- (1.2, 4); %de
\node[above] at (2.5, 4) {26};
\node[above] at (4, 4) {$V_4$}; 
\draw (3.75, 3.85) -- (.2, 2.1);  %df
\node at (3.05, 3.4)   {\small 18};
\node[above] at (1, 4) {$V_5$}; 
\draw (.83, 3.9) -- (.1, 2.1); %ef
\node[left] at (.5, 3) {14};

\node[below] at (2, -1) {{\bf a.}   $\mathcal G_S^6$};

%%%%%%%%%%%%%%%%%

\draw(11.18, 0) -- (13.9, 0); %ab 1=6, 2=3, 5-5, 4=9, 5=4, 6=8§
\node[below] at (11, 0) {$V_1$};
\node[below] at (12.5, 0) {9};
\draw(11.15, .06) -- (14.9, 2);   %ac
\node at (11.9, .5) {\small 11};
\draw (11.15, .15) -- (13.9, 3.9); %ad
\node at (12.1, 1.5)  {\small 15};
\draw (11, .2) -- (11, 3.9);  %a, e)
\node at  (11, 3) {\small 20}; 
\node[left] at (10, 2) {$V_6$}; 
\draw(10.9, .1) -- (10.1, 2);  %af
\node[left] at (10.5, 1) {14};
\draw (14.1, .1) -- (14.9, 1.9); %bc
\node[right] at  (14.5, 1) {8}; 
\node[below] at (14, 0) {$V_2$};
\draw (14, .1) -- (14, 3.85);  %bd
\node at (14, 1) {\small 12};
\draw  (13.9, .1) -- (11.1, 3.9); %be
\node at (12.1, 2.5) {\small 17};
\draw (13.82, .07) -- (10.23, 1.95);  %bf
\node at (10.56, 1.7) {\small 11 };
\draw (14.9, 2.1) -- (14.1, 3.9);  %cd
\node[right] at (14.5, 3) {14}; 
\node[right] at (15, 2){$V_3$};
\draw (14.9, 2.05) -- (11.24, 3.86); %ce
\node at (14.3, 2.35) {\small 19};
\draw (14.8, 2) -- (10.2, 2); %cf
\node at (13.3, 2) {\small 13};
\draw (13.8, 4) -- (11.2, 4); %de
\node[above] at (12.5, 4) {23};
\node[above] at (14, 4) {$V_4$}; 
\draw (13.75, 3.85) -- (10.2, 2.1);  %df
\node at (13.05, 3.45)   {\small 17};
\node[above] at (11, 4) {$V_5$}; 
\draw  (10.83, 3.9) -- (10.11, 2.1); %ef
\node[left] at (10.5, 3) {22};

\node at (7.5, 2) {=};

\node[below] at (12, -1) {{\bf b.} ${\mathcal G_{S, cpi}^6}$};

%%%%%%%%%%%%%

\draw (21.18, 0) -- (23.9, 0); %ab
\node[below] at (21, 0) {$V_1$};
\node[below] at (22.5, 0) {5};
\draw (21.15, .06) -- (24.9, 2);   %ac
\node at (21.9, .5) {\small 0};
\draw (21.15, .15) -- (23.9, 3.9); %ad
\node at (22.1, 1.5)  {0};
\draw (21, .2) -- (21, 3.9);  %a, e)
\node at  (21, 3) {\small -3}; 
\node[left] at (20, 2) {$V_6$}; 
\draw (20.9, .1) -- (20.1, 2);  %af
\node[left] at (20.5, 1) {-2};
\draw (24.1, .1) -- (24.9, 1.9); %bc
\node[right] at  (24.5, 1) {1}; 
\node[below] at (24, 0) {$V_2$};
\draw (24, .1) -- (24, 3.85);  %bd
\node at (23.92, 1) {\small -3};
\draw  (23.9, .1) -- (21.1, 3.9); %be
\node at (22.1, 2.5) {-2};
\draw (23.82, .07) -- (20.23, 1.95);  %bf
\node at (20.56, 1.7) {\small -1};
\draw  (24.9, 2.1) -- (24.1, 3.9);  %cd
\node[right] at (24.5, 3) {-1}; 
\node[right] at (25, 2){$V_3$};
\draw (24.9, 2.05) -- (21.24, 3.86); %ce
\node at (24.3, 2.35) {\small 0};
\draw  (24.8, 2) -- (20.2, 2); %cf
\node at (23.3, 2) {\small 0};
\draw  (23.8, 4) -- (21.2, 4); %de
\node[above] at (22.5, 4) {3};
\node[above] at (24, 4) {$V_4$}; 
\draw (23.75, 3.85) -- (20.2, 2.1);  %df
\node at (23.05, 3.4)   {\small 1};
\node[above] at (21, 4) {$V_5$}; 
\draw  (20.83, 3.9) -- (20.11, 2.1); %ef
\node[left] at (20.5, 3) {2};

\node at (17.5, 2) {+};

\node[below] at (22, -1) {{\bf c.} ${\mathcal G_{S, cyclic}^6}$};

 \node[below] at (12.5, -2) {{\bf Figure 9.} A $\mathcal G_S^6$};

\end{tikzpicture}

The   Fig.~9 six-alternative example   is similarly  obtained.   
The $\mathcal G^6_{S, cpi}$ weights are $\omega_1=6, \omega_2 = 3, \omega_3=5, \omega_4=9, \omega_5=14, \omega_6= 8.$ 
As required by Def.~\ref{def: S cpi},  for any rectangle in $\mathcal G^6_{S, cpi}$ (Fig.~9b),  the sums of its horizontal edges,   its vertical edges, and   its diagonals are the same.  For any five vertices, the lengths of any  $\mathcal G^6_{S, cpi}$ closed curves meeting all five vertices once are the same.  All $\mathcal G^6_{S, cpi}$  Hamiltonian paths have the same length.  Similar to Cor.~\ref{cor: cpe}, $\mathcal G^n_{S, cpi}$ can dominate the $\mathcal G^n_S$ format.

Turning to $\mathcal G^n_{S, cyclic}$, negative arc values normally are avoided  
with symmetric   costs
 because  
cycling    can generate an arbitrarily small   path length.   This problem is sidestepped here  
 because  such cycling  increases the  value that replaces  
   $T(\mathcal G^n_{S, cpi})$ in Thm.~\ref{thm: 10 done better}; e.g., if each vertex is met twice, then the value is  $2T(\mathcal G^n_{S, cpi})= 2(2\sum_{j=1}^n\omega_j)$.  As  $\mathcal G^n_{S, cyclic}$ entries indicate  
``differences from average,"   following   arcs with  negative lengths is following ``below average cost" arcs; a concept that does not exist for $\mathcal G^n_S$.  

An importance  of the reduction is that 
  the $\mathcal G^n_{S, cyclic}$ arc lengths have a distinct meaning for path lengths, so even the  GA can be successful where it would fail with $\mathcal G^n_S$.  With Fig.~8c, the GA identifies  
shortest $\mathcal G^n_{S, cyclic}$  Hamiltonian circuit  of $V_1  \stackrel{-16}{\longrightarrow} V_2  \stackrel{-1}{\longrightarrow} V_4  \stackrel{-14}{\longrightarrow} V_5  \stackrel{-4}{\longrightarrow} V_3  \stackrel{-1}{\longrightarrow} V_1$, which uses all five   negative cost arcs, has    length $-36$.   Its  $\mathcal G_S^5$ length (Thm.~\ref{thm: 10 done better})  is $T(\mathcal G_S^5) - 36 = 84-36 = 48.$ But GA is thrown off the track with $\mathcal G^5_S$ (Fig.~8a)  because of the $\mathcal G^5_{S, cpi}$ terms.

Similarly, the  GA identifies the shortest Fig.~9c  Hamiltonian path $V_1  \stackrel{-3}{\longrightarrow}  V_5  \stackrel{-2}{\longrightarrow} V_2  \stackrel{-3}{\longrightarrow} V_4  \stackrel{-1}{\longrightarrow} V_3  \stackrel{0}{\longrightarrow}  V_6  \stackrel{-2}{\longrightarrow}  V_1$  of $-11$.   
Using the $\omega_j$ values for Fig.~9b, $T(\mathcal G^6_{S, cpi})= 90,$  
 so the shortest Hamiltonian path in Fig.~9a is $11$ below this average, or $90-11=79.$  Again,  
 the GA fails for $\mathcal G^6_S$ because the $\mathcal G^6_{S, cpi}$ entries divert it.

The  construction leads to  an easily computed       
 lower bound for  Hamiltonian path lengths.

\begin{Cor}\label{cor: something} For $\mathcal G^n_S$, let the adjustment $\mathcal A(\mathcal G^n_{S, cyclic})$ be the sum of the $n$ smallest arc lengths in $\mathcal G^n_{S, cyclic}$.  All  $\mathcal G^n_S$  Hamiltonian path lengths are bounded below by $T(\mathcal G^n_S) + \mathcal A(\mathcal G^n_{S, cyclic})$.  The shortest Hamiltonian graph is bounded above by $T(\mathcal G^n_S)$, \end{Cor}

The last statement follows because $T(\mathcal G^n_S)$ is the average length of a Hamiltonian path.  Thus some Hamiltonian path length is smaller than  $T(\mathcal G^n_S)$ and $\mathcal A(\mathcal G^n_{S, cyclic})<0$ iff $\mathcal G^n_{S, cyclic}\ne 0$.   For Fig.~8,   $\mathcal A(\mathcal G^5_{S, cyclic})= -36$, so the lower bound is $84-36 =48$, which equals the length of the shortest Hamiltonian path.  With Fig.~9, $\mathcal A(\mathcal G^5_{S, cyclic}) = -12$ for the lower bound of $90-12=78$, but  the shortest Hamiltonian path has the larger length of $79.$  The reason is that the $-1$ length of $\widehat{V_2V_6}$ can not be used.  By using the four-cycle geometry, sharper estimates can be derived. 
 
Closely related to Cor.~\ref{cor: something}  is an approach to find the shortest Hamiltonian tour by  ranking $\mathcal G^n_{S, cyclic}$   arcs according to their lengths where ``smaller is better."  If marking first $n$ shortest arcs does not define a Hamiltonian tour, iteratively add arcs from this list until the marked legs do define such a closed path. (All of the shortest Hamiltonian circuits in this section were verified  in this simple manner.  This approach can be improved by using properties of the four-cycles.)   

\begin{tikzpicture}[xscale=0.5, yscale=0.5]

 \node[below] at (10.5, -2) {{\bf Figure 10.} Finding paths};
 
\node[below] at (1, 0) {$V_1$};
 
\node[below] at (4, 0) {$V_2$};

%\node at (25, 2) {$\bullet$};
\node[right] at (5, 2) {$V_3$};
%\node at (22.5, 3) {$\bullet$};
\node[above] at (2.5, 3) {$V_4$};
%\node at (20, 2) {$\bullet$};
\node[left] at (0, 2) {$V_5$};
\draw [very thick, ->] (1.1,  0) -- (3.8, 0); %ab
\draw [very thick, ->] (4.9, 1.9) -- (1.15, .15);  %ca
\draw   (1, 0) -- (2.4, 2.9); %ad
\draw  (0.9, .05) -- (0.1, 1.9); %ae
\draw  (0.15, 2) -- (3.8, .15);  %be
 \draw [very thick, <-] (2.55, 3) -- (3.9, .15);  %db
 \draw   (4, 0) -- (4.9, 1.8); %bc
\draw  (4.9, 2.1) -- (2.75, 3); %cd
\draw  [very thick, <-] (4.85, 2) -- (0.45, 2); %ce
\draw [very thick, <-] (0.2, 2.1) -- (2.3, 2.95); %ed

\node[below] at (2.4, 0) {$-{16}$};  %ab
\node[right] at (4.5, 1) {$3$}; %bc
\node[left] at (0.5, 1) {$4$}; %ae
\node[above] at (3.75, 2.5) {$2$};  %cd
\node[above] at (1.25, 2.5) {-14}; %ed

%\node[below] at (3.1, 1.3) {$0$}; %ac
\node[right] at (2.25, 1) {\small-1}; %AC
%\node[left] at (1.8, 1.75) {$2$}; %ad
\node[left] at (2.3, 1.6) {\small{13}}; %ad
%\node[below] at (1.9, 1.3) {$1$}; %eb
\node[right] at (1.4, 1) {\small 14};
%\node[right] at (3.2, 1.75) {$0$}; %bd
\node[above] at (3.3, 1.3) {\small -1}; %bd
%\node[above] at (2.5, 2) {$0$}; %ec
\node[left] at (2.9, 2.18) {\small -4}; %ec

\node[below] at (2.5, -1) {{\bf a.} From Fig.~8c}; 

%%%%%%%%%%%%%%

\draw (15.18, 0) -- (17.9, 0); %ab
\node[below] at (15, 0) {$V_1$};
\node[below] at (16.5, 0) {5};
\draw (15.15, .06) -- (18.9, 2);   %ac
\node at (15.9, .5) {\small 0};
\draw (15.15, .15) -- (17.9, 3.9); %ad
\node at (16.1, 1.5)  {0};
\draw [very thick, ->]  (15, .2) -- (15, 3.8);  %a, e)
\node at  (15, 3) {\small -3}; 
\node[left] at (14, 2) {$V_6$}; 
\draw [very thick, <-] (14.9, .1) -- (14.1, 2);  %af
\node[left] at (14.5, 1) {-2};
\draw (18.1, .1) -- (18.9, 1.9); %bc
\node[right] at  (18.5, 1) {1}; 
\node[below] at (18, 0) {$V_2$};
\draw [very thick, ->]  (18, .1) -- (18, 3.85);  %bd
\node at (17.92, 1) {\small -3};
\draw [very thick, <-]  (17.8, .2) -- (15.1, 3.8); %be
\node at (16.1, 2.5) {-2};
\draw [very thick] (17.82, .07) -- (14.23, 1.95);  %bf
\node at (14.56, 1.7) {\small -1};
\draw [very thick, <-] (18.9, 2.1) -- (18.1, 3.9);  %cd
\node[right] at (18.5, 3) {-1}; 
\node[right] at (19, 2){$V_3$};
\draw (18.9, 2.05) -- (15.24, 3.86); %ce
\node at (18.3, 2.35) {\small 0};
\draw  (18.8, 2) -- (14.2, 2); %cf
\node at (17.3, 2) {\small 0};
\draw  (17.8, 4) -- (15.2, 4); %de
\node[above] at (16.5, 4) {3};
\node[above] at (18, 4) {$V_4$}; 
\draw (17.75, 3.85) -- (14.2, 2.1);  %df
\node at (17.05, 3.4)   {\small 1};
\node[above] at (15, 4) {$V_5$}; 
\draw  (14.83, 3.9) -- (14.11, 2.1); %ef
\node[left] at (14.5, 3) {2};

\node[below] at (16, -1) {{\bf b.} From Fig.~9c.};

\end{tikzpicture}

To illustrate with Fig.~10a (from Fig.8c), just  marking the five legs with negative costs already defines the shortest Hamiltonian graph.  In Fig.~10b, the six legs with the smallest (all negative) values do not define a Hamiltonian circuit, so add an additional leg with  the next smallest  cost (here zero). The  $V_3 \stackrel{0}{\longrightarrow}  V_6$ arc completes the $V_1  \stackrel{-3}{\longrightarrow} V_5  \stackrel{-2}{\longrightarrow} V_2  \stackrel{-3}{\longrightarrow} V_4  \stackrel{-1}{\longrightarrow} V_3  \stackrel{0}{\longrightarrow} V_6  \stackrel{-2}{\longrightarrow} V_1$ circuit, which, by construction, is the shortest.

This approach applies  
to other types of closed paths. Suppose the goal in Fig.~9 is to find find the shortest  closed path that passes once through each of the four vertices $\{V_1, V_2, V_4, V_5\}$.  The six $\mathcal G^6_{S, cyclic}$  arc lengths of these vertices  are $\{-3, -3,  -2, 0, 3, 5\}$ where marking the first four on $\mathcal G^6_{S, cyclic}$   already defines the minimal closed path $V_2 \stackrel{-3}{\longrightarrow}  V_4  \stackrel{0}{\longrightarrow} V_1  \stackrel{-3}{\longrightarrow} V_5  \stackrel{-2}{\longrightarrow} V_2$ of length $-8$.   This set's $T$ value is $2(\omega_1+\omega_2+\omega_4+\omega_5),$ which, in $\mathcal G^6_{S, cpi}$,   is   the sum of the rectangle's vertical and  horizontal  legs or  64.  So the length of this shortest $\mathcal G^6_S$  closed  path over these vertices  is $64-8= 56.$

\subsection{Four cycle structure} 
A complication in determining which four-cycles define a given $\mathcal G^n_{S, cyclic}$ is that some of these four-cycles  must overlap on certain edges.  To handle this complexity, the switching, iterative approach used in the proof of Thm.~\ref{thm: s cyclic basis} is used.

\begin{Thm}\label{thm: finding S cyclic basis} To express a $\mathcal G^n_{S, cycle} \in \mathbb C^n_S$ in terms of the  basis in Eq.~\ref{eq: CS basis},  
 for $4\le j<k$, the multiple of
  $\mathbf b^n_{1, 2, j, k} \in \mathcal A^n_{1, 2}$  is $d_{j, k}$ from the $V_j  \stackrel{d_{j, k}}{\longrightarrow} V_k$ arc in $\mathcal G^n_{S, cycle}$.  (If the arc is not in the graph, its  value is zero.)  After determining the multiple of a basis vector, remove the associated four-cycle from the graph.   
 In what  remains, the multiple of the top $\mathbf b^n_{1, 2, n, 3} \in \mathcal B^n_{1, 2}$ is the negative of the $d_{2, n}$ value in the of $V_2 \stackrel{d_{2, n}}{\longrightarrow} V_n$ arc in the   reduced graph.\footnote{The associated arc for $\mathbf b^n_{1, 2, n, 3}$  is $V_1 \stackrel{1}{\longrightarrow} V_2 \stackrel{-1}{\longrightarrow} V_n \stackrel{1}{\longrightarrow} V_3 \stackrel{-1}{\longrightarrow} V_1$, so for $V_2 \stackrel{d_{2, n}}{\longrightarrow} V_n$ to hold, the   coefficient for  $\mathbf b^n_{1, 2, n, 3}$ must be the negative of $d_{2, n}$.}  After removing this four-cycle, the top remaining $\mathcal A^n_{1, 2}$ vector is   $\mathbf b^n_{1, 2, 3, n}$; its coefficient is the length in $V_3 \stackrel{d_{3, n}}{\longrightarrow} V_n$ in the reduced graph, which leaves $n-1$ as  the largest remaining index in the reduced graph.   In general, if the largest remaining index is s,  the multiple of the  top remaining $\mathcal B^n_{1, 2}$ vector, $\mathbf b^n_{1, 2, s, 3}$, is the negative of    $d_{2, s}$   from the reduced graph's $V_2 \stackrel{d_{2, s}}{\longrightarrow} V_s$ arc.  The top of the remaining $\mathcal A^n_{1, 2}$ vectors is $\mathbf b^n_{1, 2, 3, s}$; its multiple is the $d_{3, s}$ value of the $V_3 \stackrel{d_{3, s}}{\longrightarrow} V_s$ arc in the reduced graph. \end{Thm} \smallskip
 
{\em Proof,}  The proof is essentially  
    that of Thm.~\ref{thm: s cyclic basis};     removing basis vectors in the specified manner leaves, at each stage, a single $d_{u, v}$ value of a certain type.  Because $\mathbb C^n_S$  is the sum of these four-cycles, the existence of  this $d_{u, v}\ne 0$ requires the associated $\mathbf b^n_{1, 2 , k, s}$ to be in the decomposition; the form of this four-cycle requires  $d_{u, v}$  to be the vector's multiple.  A difference is that if $d_{u, v}$ identifies a vector from $\mathcal B^n_{1, 2}$, the multiple is the negative of $d_{u, v}$, as required by the form of the associated four-cycle.  If the vector is from $\mathcal A^n_{1, 2}$, then $d_{u, v}$ is the multiple.  $\square$  \smallskip  
    
    Using this approach, the four cycles of $\mathcal G^5_{S, cyclic}$ in Fig.~8c are $-14\mathbf b^5_{1, 2, {\mathbf 4, 5}}$, $-\mathbf b^5_{1, \mathbf{2, 5,} 3}$, $10\mathbf b^5_{1, 2, \mathbf{3, 5}}$, $-13\mathbf b^5_{1, \mathbf{2, 3,} 4}$ and $15\mathbf b^5_{1, 2, 4, 3}$.

\subsection{Extensions}  With minor   modifications, all other results developed in Sect.~2  for the asymmetric $\mathbb G^n_A$ transfer to the symmetric $\mathbb G^n_S$.   For instance, to analyze connected and closed path properties that involve a subset of vertices,  carry out the above with that subset.  Other samples follow.

\begin{Thm} \label{thm: s gen paths}  Consider the class of paths starting at $V_j$ and ending at $V_k$ that pass through vertices $\{V_i\}_{i\in \mathcal D}$ where, for each $i\in \mathcal D$, the path passes through $V_i$  $\kappa_i$ times.  The length of such a path in $\mathcal G^n_S$ is its path length in $\mathcal G^n_{S, cyclic}$ plus $(\omega_j+\omega_k) + 2 \sum_{i\in \mathcal D} \kappa_i\omega_i$.  \end{Thm}

As an example, consider all paths in Fig.~9 that start at $V_1$, end in $V_5$ and pass through each of $V_2, V_3, V_4$ twice.  According to the weights of $\mathcal G^6_{S, cpi}$, the length of any of these paths in $\mathcal G^6_S$ is its length in $\mathcal G^6_{S, cyclic}$ plus its $\mathcal G^6_{S, cpi}$ length of $(6+ 4) + 4(3+5+4).$

 Graphs with incomplete symmetric costs  are handled  the same way as  
 in  Sect.~2.  That is, complete the graph by adding arcs of any desired length to obtain $\tilde{\mathcal G}_S^n$. 
   Then compute $\tilde{\mathcal G}^n_{S, cpi}$ and  $\tilde{\mathcal G}^n_{S, cyclic}$   For incomplete graphs,     $\infty$ often is assigned to  inadmissible arcs;   
  do so only with   $\tilde{\mathcal G}^n_{S, cyclic}.$

\begin{Thm}\label{thm: incomplete symmetric}  For an incomplete symmetric graph $\mathcal G^n_S$,  let $\tilde{\mathcal G_S^n}$ include  
the missing $\mathcal G_S^n$  arcs where each has    
an arbitrary selected  length.  Compute $T(\tilde{\mathcal G_S^n})$ and $\tilde{\mathcal G}^n_{S, cyclic}.$  The length of a $\mathcal G_S^n$  Hamiltonian path  
  is  $T(\tilde{\mathcal G}_S^n)$ plus its  $\tilde{\mathcal G}^n_{S, cyclic}$ length. \end{Thm}

  Computations can be simplified by adding arcs of zero length so that   the   $S_S(V_j)$ values for $\mathcal G_S^n$ and $\tilde{\mathcal G_S^n}$   agree, and  
 $T(\mathcal G_S^n) = T(\tilde{\mathcal G_S^n}).$   
   \smallskip

\noindent{\em Proof:}   A  Hamiltonian  path length in $\mathcal G_S^n$ is the same   in $\tilde{\mathcal G_S^n}$, which  equals   $T(\tilde{\mathcal G_S^n})$ plus its length in $\tilde{\mathcal G^n}_{S, cyclic}.$  The result follows.  $\square$

 As $\widehat{V_1V_4}$ and $\widehat{V_2V_5}$ are not admitted in Fig.~11,    include them in Fig.~11a    with zero lengths (the two dashed arcs in Fig.~11a).   
 The $S_S(V_j)$ values of   $\tilde{\mathcal G_S^n}$ are $S_S(V_1) = 27, S_S(V_2) = 27, S_S(V_3) = 39, S_S(V_4) = 23, S_S(V_5) =  23, S_S(V_6) = 51$. Thus  
  $T(\tilde{\mathcal G_S^n}) = 38,$   
$\omega_1=\frac14[27-19]=2, \omega_2=2, \omega_3= 5, \omega_4= 1, \omega_5=1, \omega_6 = 8, $ and Figs.~10b, c follow.   The two inadmissible  Fig.~11c arcs  (both with length   of $-3$)  could be dropped or, as in Fig.~11c,   replaced     with $\infty.$   The shortest  $\mathcal G^6_{S, cyclic}$ Hamiltonian  path $V_1 \stackrel{-2}{\longrightarrow} V_3 \stackrel{0}{\longrightarrow} V_2 \stackrel{-1}{\longrightarrow} V_6  \stackrel{-2}{\longrightarrow} V_4 \stackrel{2}{\longrightarrow} V_5 \stackrel{0}{\longrightarrow} V_1$ of length $-3$, which  includes all allowed  arcs with negative  
   costs, can be found in the above described manner.   In $\mathcal G_S^6$ this path has the ``below average"  length of $T(\mathcal G_S^6) - 3 =  
   35.$

 \begin{tikzpicture}[xscale=0.5, yscale=0.52]

\draw (1.18, 0) -- (3.9, 0); %ab
\node[below] at (1, 0) {$V_1$};
\node[below] at (2.5, 0) {7};
\draw (1.15, .06) -- (4.9, 2);   %ac
\node at (1.9, .5) {\small 5};
\draw[dashed] (1.15, .15) -- (3.9, 3.9); %ad
\node at (2.1, 1.5)  {\small 0};
\draw (1, .2) -- (1, 3.9);  %a, e)
\node at  (1, 3) {\small 3}; 
\node[left] at (0, 2) {$V_6$}; 
\draw (.9, .1) -- (.1, 2);  %af
\node[left] at (.5, 1) {12};
\draw (4.1, .1) -- (4.9, 1.9); %bc
\node[right] at  (4.5, 1) {7}; 
\node[below] at (4, 0) {$V_2$};
\draw (4, .1) -- (4, 3.85);  %bd
\node at (4, 1) {\small 4};
\draw[dashed]  (3.9, .1) -- (1.1, 3.9); %be
\node at (2.1, 2.5) {\small 0};
\draw (3.82, .07) -- (.23, 1.95);  %bf
\node at (.56, 1.7) {\small 9};
\draw (4.9, 2.1) -- (4.1, 3.9);  %cd
\node[right] at (4.5, 3) {8}; 
\node[right] at (5, 2){$V_3$};
\draw (4.9, 2.05) -- (1.24, 3.86); %ce
\node at (4.3, 2.35) {\small 6};
\draw (4.8, 2) -- (.2, 2); %cf
\node at (3.3, 2) {\small 13};
\draw (3.8, 4) -- (1.2, 4); %de
\node[above] at (2.5, 4) {4};
\node[above] at (4, 4) {$V_4$}; 
\draw (3.75, 3.85) -- (.2, 2.1);  %df
\node at (3.05, 3.4)   {\small 7};
\node[above] at (1, 4) {$V_5$}; 
\draw (.83, 3.9) -- (.1, 2.1); %ef
\node[left] at (.5, 3) {10};

\node[below] at (2, -1) {{\bf a.} From $\mathcal G_S^6$ to $\widetilde{\mathcal G_S^6}$};

%%%%%%%%%%%%%%%%%

\draw(11.18, 0) -- (13.9, 0); %ab
\node[below] at (11, 0) {$V_1$};
\node[below] at (12.5, 0) {4};
\draw(11.15, .06) -- (14.9, 2);   %ac
\node at (11.9, .5) {\small 7};
\draw[dashed] (11.15, .15) -- (13.9, 3.9); %ad
\node at (12.1, 1.5)  {\small 3};
\draw (11, .2) -- (11, 3.9);  %a, e)
\node at  (11, 3) {\small 3}; 
\node[left] at (10, 2) {$V_6$}; 
\draw(10.9, .1) -- (10.1, 2);  %af
\node[left] at (10.5, 1) {10};
\draw (14.1, .1) -- (14.9, 1.9); %bc
\node[right] at  (14.5, 1) {7}; 
\node[below] at (14, 0) {$V_2$};
\draw (14, .1) -- (14, 3.85);  %bd
\node at (14, 1) {\small 3};
\draw[dashed]  (13.9, .1) -- (11.1, 3.9); %be
\node at (12.1, 2.5) {\small 3};
\draw (13.82, .07) -- (10.23, 1.95);  %bf
\node at (10.56, 1.7) {\small 10 };
\draw (14.9, 2.1) -- (14.1, 3.9);  %cd
\node[right] at (14.5, 3) {6}; 
\node[right] at (15, 2){$V_3$};
\draw (14.9, 2.05) -- (11.24, 3.86); %ce
\node at (14.3, 2.35) {\small 6};
\draw (14.8, 2) -- (10.2, 2); %cf
\node at (13.3, 2) {\small 13};
\draw (13.8, 4) -- (11.2, 4); %de
\node[above] at (12.5, 4) {2};
\node[above] at (14, 4) {$V_4$}; 
\draw (13.75, 3.85) -- (10.2, 2.1);  %df
\node at (13.05, 3.45)   {\small 9};
\node[above] at (11, 4) {$V_5$}; 
\draw  (10.83, 3.9) -- (10.11, 2.1); %ef
\node[left] at (10.5, 3) {9};

\node at (7.5, 2) {=};

\node[below] at (12, -1) {{\bf b.} $\widetilde{\mathcal G_{S, cpi}^6}$};

%%%%%%%%%%%%%

\draw (21.18, 0) -- (23.9, 0); %ab
\node[below] at (21, 0) {$V_1$};
\node[below] at (22.5, 0) {3};
\draw (21.15, .06) -- (24.9, 2);   %ac
\node at (21.9, .5) {\small -2};
\draw[dashed] (21.15, .15) -- (23.9, 3.9); %ad
\node at (22.1, 1.5)  {$\infty$};
\draw (21, .2) -- (21, 3.9);  %a, e)
\node at  (21, 3) {\small 0}; 
\node[left] at (20, 2) {$V_6$}; 
\draw (20.9, .1) -- (20.1, 2);  %af
\node[left] at (20.5, 1) {2};
\draw (24.1, .1) -- (24.9, 1.9); %bc
\node[right] at  (24.5, 1) {0}; 
\node[below] at (24, 0) {$V_2$};
\draw (24, .1) -- (24, 3.85);  %bd
\node at (23.92, 1) {\small 1};
\draw[dashed]  (23.9, .1) -- (21.1, 3.9); %be
\node at (22.1, 2.5) {$\infty$};
\draw (23.82, .07) -- (20.23, 1.95);  %bf
\node at (20.56, 1.7) {\small -1};
\draw  (24.9, 2.1) -- (24.1, 3.9);  %cd
\node[right] at (24.5, 3) {2}; 
\node[right] at (25, 2){$V_3$};
\draw (24.9, 2.05) -- (21.24, 3.86); %ce
\node at (24.3, 2.35) {\small 0};
\draw  (24.8, 2) -- (20.2, 2); %cf
\node at (23.3, 2) {\small 0};
\draw  (23.8, 4) -- (21.2, 4); %de
\node[above] at (22.5, 4) {2};
\node[above] at (24, 4) {$V_4$}; 
\draw (23.75, 3.85) -- (20.2, 2.1);  %df
\node at (23.05, 3.4)   {\small  -2};
\node[above] at (21, 4) {$V_5$}; 
\draw  (20.83, 3.9) -- (20.11, 2.1); %ef
\node[left] at (20.5, 3) {1};

\node at (17.5, 2) {+};

\node[below] at (22, -1) {{\bf c.} $\widetilde{\mathcal G_{S, cyclic}^6}$};

 \node[below] at (12.5, -2.2) {{\bf Figure 11.} An incomplete $\mathcal G_S^6$};

\end{tikzpicture}

While the method associated with Fig.~10 is more efficient than the GA, it is worth using the structure of the decomposition to explain certain  GA traits.  Parallel to the $\mathcal G^n_{A, cyclic}$ concern whether a vertex can be a source, a  $\mathcal G^n_{S, cyclic}$ worry is whether all of a vertex's  legs are positive.   But  $S_S(V_j)=0$,  so this cannot happen.   Thus, if  all  $\mathcal G^n_{S, cyclic}$ vertices have an arc with non-zero length, then the number of negative length arcs is bounded below by $\frac n2$ and above by $\frac{n(n-2)}2$.  Dropping  $\mathcal G^n_{S, cpi}$ eliminates one   GA difficulty, but another   
  is caused by the number of options.   
 This can be seen with Fig.~12a, which is given by $\{-x_1\mathbf b^6_{1, 2, 5, 6} -x_2\mathbf b^6_{1, 2, 6, 5}\}$, $\{ - y_1 \mathbf b^6_{1, 2, 3, 4}, -y_2\mathbf b^6_{1, 2, 4, 3}\}$, and $\{-z_1\mathbf b^6_{3, 4, 5, 6}, -z_2\mathbf b^6_{3, 4, 6, 5}\},$ where each bracket defines a rectangle.   Thus $u=-(x_1+x_2) -(y_1+y_2)$, $v=-(y_1+y_2)-(z_1+z_2)$, and $w=-(x_1+x_2)-(z_1+z_2)$.  Should the $x$'s,   $y$'s, and   $z$'s have positive values, the graph has $\frac n2 =3$  negative  and  $\frac{n(n-2)}2=12$ positive arc lengths.
 
 \begin{tikzpicture}[xscale=0.45, yscale=0.52]

\draw (1.18, 0) -- (3.9, 0); %ab
\node[below] at (1, 0) {$V_1$};
\node[below] at (2.5, 0) {u};
\draw (1.15, .06) -- (4.9, 2);   %ac
\node at (1.9, .5) {$y_2$};
\draw (1.15, .15) -- (3.9, 3.9); %ad
\node at (2.1, 1.5)  {$y_1$};
\draw (1, .2) -- (1, 3.9);  %a, e)
\node at  (1, 3) {$x_2$}; 
\node[left] at (0, 2) {$V_6$}; 
\draw (.9, .1) -- (.1, 2);  %af
\node[left] at (.5, 1) {$x_1$};
\draw (4.1, .1) -- (4.9, 1.9); %bc
\node[right] at  (4.5, 1) {$y_1$}; 
\node[below] at (4, 0) {$V_2$};
\draw (4, .1) -- (4, 3.85);  %bd
\node at (4, 1) {$y_2$};
\draw  (3.9, .1) -- (1.1, 3.9); %be
\node at (2.1, 2.5) {$x_1$};
\draw (3.82, .07) -- (.23, 1.95);  %bf
\node at (.56, 1.7) {$x_2$};
\draw (4.9, 2.1) -- (4.1, 3.9);  %cd
\node[right] at (4.5, 3) {v}; 
\node[right] at (5, 2){$V_3$};
\draw (4.9, 2.05) -- (1.24, 3.86); %ce
\node at (4.3, 2.35) {$z_2$};
\draw (4.7, 2) -- (.2, 2); %cf
\node at (3.3, 2) {$z_1$};
\draw (3.8, 4) -- (1.2, 4); %de
\node[above] at (2.5, 4) {$z_1$};
\node[above] at (4, 4) {$V_4$}; 
\draw (3.75, 3.85) -- (.2, 2.1);  %df
\node at (3.05, 3.4)   {$z_2$};
\node[above] at (1, 4) {$V_5$}; 
\draw (.83, 3.9) -- (.1, 2.1); %ef
\node[left] at (.5, 3) {w};

\node[below] at (2, -1) {{\bf a.} A choice of $\mathcal G^6_{S, cyclic}$};

%%%%%%%%%%%%%%%%%

\draw[thick, ->] (9.18, 0) -- (11.9, 0); %ab 1=6, 2=3, 5-5, 4=9, 5=4, 6=8§
\node[below] at (9, 0) {$V_1$};
\node[below] at (10.5, 0) {u};
\draw[dashed](9.15, .06) -- (12.9, 2);   %ac
\node at (9.9, .5) {$y_2$};
\draw[thick, <-] (9.15, .15) -- (11.9, 3.9); %ad
\node at (10.1, 1.5)  {$y_1$};
\draw[dashed] (9, .2) -- (9, 3.9);  %a, e)
\node at  (9, 3) {$x_2$}; 
\node[left] at (8, 2) {$V_6$}; 
\draw[dashed] (8.9, .1) -- (8.1, 2);  %af
\node[left] at (8.5, 1) {$x_1$};
\draw[dashed] (12.1, .1) -- (12.9, 1.9); %bc
\node[right] at  (12.5, 1) {$y_1$}; 
\node[below] at (12, 0) {$V_2$};
\draw[dashed] (12, .1) -- (12, 3.85);  %bd
\node at (12, 1) {$y_2$};
\draw[thick, ->]  (11.9, .1) -- (9.1, 3.9); %be
\node at (10.1, 2.5) {$x_1$};
\draw[dashed] (11.82, .07) -- (8.23, 1.95);  %bf
\node at (8.56, 1.7) {$x_2$ };
\draw[thick, ->] (12.9, 2.1) -- (12.1, 3.9);  %cd
\node[right] at (12.5, 3) {v}; 
\node[right] at (13, 2){$V_3$};
\draw[dashed, ->] (12.9, 2.05) -- (9.24, 3.86); %ce
\node at (12.3, 2.35) {$z_2$};
\draw[thick, <-] (12.8, 2) -- (8.2, 2); %cf
\node at (11.3, 2) {$z_1$};
\draw[dashed] (11.8, 4) -- (9.2, 4); %de
\node[above] at (10.5, 4) {$z_1$};
\node[above] at (12, 4) {$V_4$}; 
 
\draw[dashed, ->] (11.75, 3.85) -- (8.2, 2.1);  %df
\node at (11.05, 3.45)   {$z_2$};
\node[above] at (9, 4) {$V_5$}; 
\node at (9, 4) {$\bullet$}; 
\draw[thick, ->]  (8.83, 3.9) -- (8.11, 2.1); %ef
\node[left] at (8.5, 3) {w};

\node[below] at (10, -1) {{\bf b.} Selecting $z_1$};

%%%%%%%%%%%%%

\draw[thick, ->] (17.18, 0) -- (19.9, 0); %ab
\node[below] at (17, 0) {$V_1$};
\node[below] at (18.5, 0) {u};
\draw[thick, <-] (17.15, .06) -- (20.9, 2);   %ac
\node at (17.9, .5) {$y_2$};
\draw[dashed] (17.15, .15) -- (19.9, 3.9); %ad
\node at (18.1, 1.5)  {$y_1$};
\draw[dashed] (17, .2) -- (17, 3.9);  %a, e)
\node at  (17, 3) {$x_2$}; 
\node[left] at (16, 2) {$V_6$}; 
\draw[dashed] (16.9, .1) -- (16.1, 2);  %af
\node[left] at (16.5, 1) {$x_1$};
\draw[dashed] (20.1, .1) -- (20.9, 1.9); %bc
\node[right] at  (20.5, 1) {$y_1$}; 
\node[below] at (20, 0) {$V_2$};
\draw[dashed] (20, .1) -- (20, 3.85);  %bd
\node at (19.92, 1) {$y_2$};
\draw[thick, ->]  (19.9, .1) -- (17.1, 3.9); %be
\node at (18.1, 2.5) {$x_1$};
\draw[dashed] (19.82, .07) -- (16.23, 1.95);  %bf
\node at (16.56, 1.7) {$x_2$};
\draw[thick, <-]  (20.9, 2.1) -- (20.1, 3.9);  %cd
\node[right] at (20.5, 3) {v}; 
\node[right] at (21, 2){$V_3$};
\draw[dashed] (20.8, 2.1) -- (17.24, 3.86); %ce
\node at (20.3, 2.35) {$z_2$};
\draw[dashed]  (20.8, 2) -- (16.3, 2); %cf
\node at (19.3, 2) {$z_1$};
\draw[dashed]  (19.8, 4) -- (17.2, 4); %de
\node[above] at (18.5, 4) {$z_1$};
\node[above] at (20, 4) {$V_4$}; 
\draw[thick, <-] (19.75, 3.85) -- (16.2, 2.1);  %df
\node at (19.05, 3.4)   {$z_2$};
\node[above] at (17, 4) {$V_5$}; 
\draw[thick, ->]  (16.83, 3.9) -- (16.11, 2.1); %ef
\node[left] at (16.5, 3) {w};

\node[below] at (18, -1) {{\bf c.} Selecting $z_2$};

%%%%%%%%%%%%%%%%%%%%%%

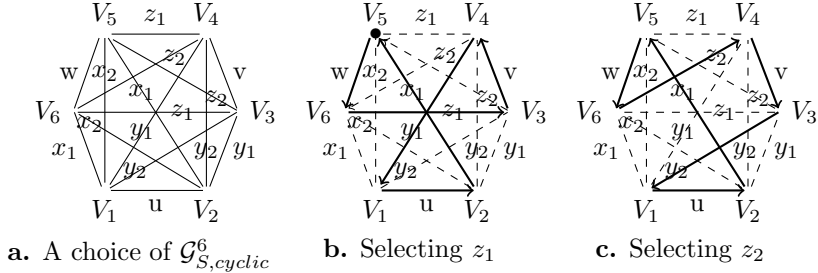
 \node[below] at (10, -2.2) {{\bf Figure 12.} Potential failures of the Greedy Algorithm};

\end{tikzpicture}

 It follows from the material following Cor.~\ref{cor: something} that a shortest Hamiltonian path must include the three arcs with negative lengths (of u, v, and w).  Thus the four-cycle symmetry requires all Hamiltonian paths  of Fig.~12a with  these negative length arcs
  to have one of only four sizes where $u+v+w$ is supplemented by $x_1+y_1+z_1, x_1+y_2+z_2, x_2 + y_1+z_1$ or $x_2+y_2+z_2$.   Starting the GA at $V_1$ and assuming that $x_1$ is smaller than $x_2, y_1, y_2$, the first three moves (Figs.~12b, c) are
 $V_1 \stackrel{u}{\longrightarrow} V_2 \stackrel{x_{1}}{\longrightarrow} V_5 \stackrel{w}{\longrightarrow}  V_6$.  The next move is to the $\widehat{V_3V_4}$ arc.  Whichever way $\widehat{V_3V_4}$ is entered determines the last arc.

Should $z_1<z_2$, then Fig.~12b represents the fourth GA step  $V_6 \stackrel{z_1}{\longrightarrow} V_3$.   But should $y_1$ be much larger than $y_2$,   the Fig.~12c route would be shorter.  As the arrangement of the  
 $\mathcal G^n_{S cyclic}$   four-cycles  
  can  affect the success of  an algorithm,  
  the algebra of these four-cycles needs to be better understood.

  For small values of $n$, and theoretically for all values,  the basis for $\mathbb C^n_S$ exhibits all  possible   $\mathcal G^n_{S, cyclic}$ choices.      To illustrate an application,  
    if arc costs represent  Euclidean distances in the planar problem and the triangle inequality is satisfied, then  minimal Hamiltonian paths cannot have a self intersection.  It is reasonable to wonder  where else  
   does the triangle inequality ensure this behavior. 
    Here, the structure of $\mathcal G^n_{S, cpi}$ plays a role.

  \begin{Thm}\label{thm: triangle}   For $n\ge 3$, $\mathcal G^n_{S, cpi}$ satisfies the triangle inequality iff all weights $\omega_j$ are non-negative.    \end{Thm}
  
     According to Eq.~\ref{def: w eq},   $\omega_k<0$ represents  
 where the associated $S_S(V_k)$ is bounded above by $\frac12T(\mathcal G^n_S)$.  That is,  the average  of  arc lengths attached to $V_k$ is much smaller than average over the graph.\smallskip
  
 {\em Proof:}  The $\mathcal G^n_{S, cpi}$ arc length for $\widehat{V_iV_j}$ is $\omega_i + \omega_j$.  For a triplet, the length of the two arcs $\widehat{V_iV_j}$ and $\widehat{V_jV_k}$ is $\omega_i +2\omega_j +\omega_k$, which differs from the $\widehat{V_iV_k}$ arc length by $2\omega_j$.  Thus, the triangle inequality is satisfied iff $\omega_j\ge0$.  This must hold for all legs of  all triplets, so  $\omega_j\ge0$ for all  $j$. $\square$

Turning to $n=4$, Fig.~13a represents all possible $\mathcal G^4_{S, cyclic}$ structures (with $u\mathbf b^4_{1, 4, 2, 3}+v\mathbf b^4_{1, 2, 3, 4}$), so it characterizes all  closed path properties and their lengths for   $\mathcal G^4_S\in \mathbb G^4_s$.    Assuming the rectangular Fig.~13a faithfully represents the geometry of a considered concern (e.g., using actual rather than Euclidean costs), the issue is to understand which $(u, v)$ values (that is, which  $\mathcal G^4_{S, cyclic}$)   require a shortest Hamiltonian path to avoid    the crossing diagonals.  
 The answer follows  from Fig.~13a as it requires the length of the two diagonals to be greater than that of the two vertical and the two horizontal edges, or $  -u>v, \,  -u> u-v$; this is the open, unbounded,  shaded Fig.~13b region.
 
 \begin{tikzpicture}[xscale=0.5, yscale=0.5]

\draw (-24, 0) -- (-21, 0);
\draw (-24, 2.6) -- (-21, 2.6);
\draw (-24, 0.1) -- (-24, 2.5);
\draw (-21, 0.1) -- (-21, 2.5); 
\draw (-23.9, 0.1) -- (-21.1, 2.5);
\draw (-23.9, 2.5) -- (-21.1, .1);
\node[above] at (-22.5, 2.6) {v};
\node[below] at (-22.5, 0) {v};
\node[left] at (-24, 1.3) {u-v};
\node[right] at (-21, 1.3) {u-v};
\node at (-22, 2) {-u};
\node at (-23.2, 2) {-u};
\node[left] at (-24, 0) {$V_1$};
\node[left] at (-24, 2.5) {$V_4$};
\node[right] at (-21, 0) {$V_2$}; 
\node[right] at (-21, 2.5) {$V_3$}; 
\node[below] at (-21.85, -1) {{\bf a.} All $\mathcal G^4_{S, cyclic}$};  % $u \mathbf b^4_{1, 4, 2, 3} +v\mathbf b^4_{1, 2, 3, 4}$};

%%%%%%%%%%%

\draw[thick]  (-18.5, 2) -- (-13.5, 2);
\draw[thick]  (-16, 0) -- (-16, 3.5);  %-16, 2
\draw (-18, 4) -- (-14, 0);  %v=-u
\draw (-15, 4) --  (-17, 0);         %v=2u
\node[right] at (-15, 1) {$v=-u$};
\node[right] at (-15.7, 3) {$v=2u$};
\draw  [gray, fill=gray] (-17, 0) -- (-16, 2) -- (-17, 3) -- (-17, 0);
\node[below] at (-15.1, -1) {{\bf b.} No diagonal crossings}; 
\node[right] at  (-13.5, 1.9) {u axis};
\node[above] at (-16, 3.5) {v axis};

%%%%%%%%%%%%%%

\draw[thick] (-10, 2) -- (-5, 2);
\draw[thick] (-7.5, -.3) -- (-7.5, 3.5);   %-7.5, 2

\draw (-8.5, -.5) -- (-8.5, 4);
\node[above] at (-8.5, 3.5) {$u=-\omega$}; 
\draw (-9.7, 3) -- (-4.6, 3);
\node [left] at (-9.7, 3) {$v=\omega$};

\draw (-8.8, -.3) -- (-5, 3.4);   %v=u-w
\node[above] at (-5, 3.12) {$v=u-\omega$}; 
\draw  [gray, fill=gray] (-8.5, 0) -- (-5.48, 3) -- (-8.5, 3)  -- (-8.5, 0);  %(-7.5, 1) -- (-7.5, 3) -- 
\node[below]  at (-7.5, -1) {{\bf c.} Triangle inequality};
\node[right] at (-5, 1.9) {u axis};
\node[below] at (-7.5, 0) {v axis};

%%%%%%%%%%%%

\node[below] at (0, 0) {$V_1$};
%%\node at (24, 0) {$\bullet$};
\node[below] at (3, 0) {$V_2$};

\node[right] at (4, 2) {$V_3$};
 
\node[above] at (1.5, 3) {$V_4$};
 
\node[left] at (-1, 2) {$V_5$};
\draw (.1,  0) -- (2.8, 0); %ab
\draw (3.9, 1.9) -- (0.15, .15);  %ca
\draw  (0, 0) -- (1.4, 2.9); %ad
\draw  (-.1, .05) -- (-.9, 1.9); %ae
\draw  (-.85, 2) -- (2.8, .15);  %be
 \draw  (1.55, 3) -- (2.9, .15);  %db
 \draw   (3, 0) -- (3.9, 1.8); %bc
\draw  (3.9, 2.1) -- (1.75, 3); %cd
\draw  (3.85, 2) -- (-.55, 2); %ce
\draw (-.8, 2.1) -- (1.3, 2.95); %ed

\node[below] at (1.4, 0) {\small{v+y-z}};  %ab
\node[right] at (2.7, 1) {\small{-v+w+z}}; %bc
\node[left] at (-.3, 1) {\small{x-y+z}}; %ae
\node[above] at (2.95, 2.3) {\small{v-w+x}};  %cd
\node[above] at (0.05, 2.3) {\small{w-x+y}}; %ed

%\node[below] at (3.1, 1.3) {$0$}; %ac
\node[right] at (1.54, .99) {\small -x}; %AC
%\node[left] at (1.8, 1.75) {$2$}; %ad
\node[left] at (1.3, 1.75) {\small -v}; %ad
%\node[below] at (1.9, 1.3) {$1$}; %eb
\node[right] at (0.56, .99) {\small -w};
%\node[right] at (3.2, 1.75) {$0$}; %bd
\node[above] at (2.23, 1.3) {\small -y}; %bd
%\node[above] at (2.5, 2) {$0$}; %ec
\node[left] at (1.95, 2.1) {\small -z}; %ec

\node[below] at (1.7, -1) {{\bf d.} All   $\mathcal G^5_{S, cyclic}$ choices}; 
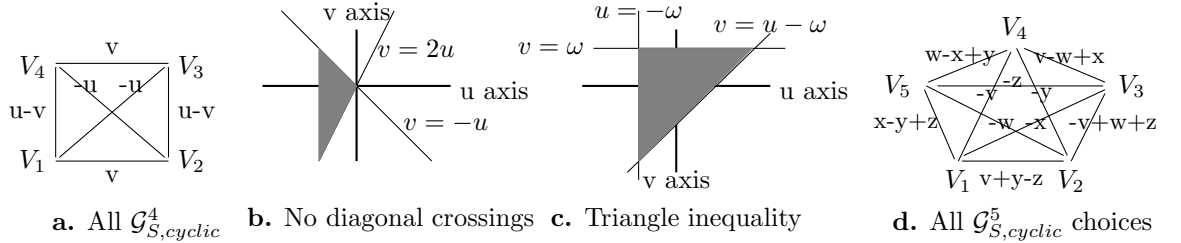
\node[below] at (-8.5, -2.2) {{\bf Figure 13.} Finding properties of $\mathcal G^n_{S, cyclic}$};

\end{tikzpicture}

To compare this wedge with what happens should $\mathcal G^4_S$ satisfy the triangle inequality, it follows from Thm.~\ref{thm: triangle} that  
  in a $\{V_i, V_j, V_k\}$ triangle, the sum of the $\widehat{V_iV_k}$
  and $\widehat{V_kV_j}$ leg lengths in $\mathcal G^4_{S, cyclic}$  plus $2\omega_k$ must be an upper bound for  $\widehat{V_iV_j}$'s $\mathcal G^4_{S, cyclic}$ leg length.  Applying this to the three triangles where vertex  $V_k$ is off the triangle's  compared edge  leads to   \begin{equation}\label{eq: for triangle} \omega_k\ge -u, \quad \omega_k \ge v, \quad \omega_k \ge u-v.\end{equation}  
   
 \begin{Thm} For $\mathcal G^4_S$, if any $\omega_j<0$, then $\mathcal G^4_S$ does not satisfy the triangle inequality.  Let $\omega=\min(\omega_1, \omega_2, \omega_3, \omega_4)$.  The region where $\mathcal G^4_S$  satisfies the triangle inequality is defined  by substituting $\omega$ for $\omega_k$ in Eq.~\ref{eq: for triangle}; it   is depicted by the closed shaded triangle in Fig.~13c. \end{Thm}

 {\em Proof:} If $\omega_k<0$. then Eq.~\ref{eq: for triangle} cannot be satisfied.  The remainder follows from the above. $\square$\smallskip
 
 If $\omega=0$,  the triangle inequality is satisfied only for  $u=v=0$, which is $\mathcal G^4_{S, cyclic}=0$ so $\mathcal G^4_S=\mathcal G^4_{S, cpi}$.  Both the triangle inequality and the non-crossing of the diagonals in the shortest Hamiltonian  path   hold  in the  
  intersection of the shaded portions of Figs.~13b, c; this  is the Fig.~13b shaded triangle limited on the left by $u\ge-\omega$.  What remains are  
   regions (i.e., choices of $\mathcal G^4_{S, cyclic}$) where the triangle inequality is satisfied but the shortest Hamiltonian path includes the diagonals, and  
   a sizable region (the shaded Fig.~13b region for  $u<-\omega$)  where  the diagonals are not in the shortest Hamiltonian circuit and the triangle inequality is not satisfied.

  Results for $\mathcal G^5_{S, cyclic}$ follow in a similar manner.    The basis for Fig.~13d, which captures all $\mathcal G^5_{S, cyclic}$  behaviors,  is  
  $\{v\mathbf b^5_{1, 2, 3, 4} + w\mathbf b^5_{2, 3, 4, 5} + x\mathbf b^5_{3, 4, 5, 1} +$ $ y\mathbf b^5_{4, 5, 1, 2} + z\mathbf b^5_{5, 1, 2, 3}\}$.

  \section{Graphs with general asymmetric costs}\label{sect: pareto}
  
 Other  systems  can be similarly reduced.   
 Graphs where  all  
 closed paths have a fixed  
 length    identify   components of $\mathcal G^n$ entries   that frustrate  
  a closed path analysis.    The subspace's  normal bundle measures  
 deviations   from    neutrality,  so  it  is   critical  when  determining closed path properties. 
  
  None of this is necessary for the standard space of graphs with  asymmetric costs.   
  The reason is that, for each pair,  the arc lengths $V_j \stackrel{x}{\longrightarrow} V_k$ and $V_k \stackrel{y}{\longrightarrow} V_j$    can be represented as  an  \{average cost, excess cost\} pair; e.g.,  $\{a=\frac{x+y}2,  V_j \stackrel{x-a}{\longrightarrow} V_k\}$.  By    applying the   Sects.~2 and 3 approaches to each component,   the above results about incomplete graphs, path lengths, etc.,  transfer.

  \begin{tikzpicture}[xscale=0.45, yscale=0.55]

\node[below] at (1, -.10) {$V_1$};
 
\node[below] at (4, -.10) {$V_2$};

\node[right] at (5, 2) {$V_3$};
%\node at (22.5, 3) {$\bullet$};
\node[above] at (2.5, 3) {$V_4$};
%\node at (20, 2) {$\bullet$};
\node[left] at (0, 2) {$V_5$};
\draw[<-] (1.1,  0) -- (3.8, 0); %ab
\draw[->] (1.1, -.2) -- (3.8, -.2);
\draw [->] (4.9, 1.9) -- (1.25, .15);  %ca
\draw [<-] (4.8, 1.75) -- (1.25, .05);
\draw [->]  (.9, .1) -- (2.25, 2.8); %ad
\draw [<-]  (1.19, .19) -- (2.4, 2.8); 
\draw [<-] (.9, .05) -- (.1, 1.9); %ae
\draw [->] (.7, .05) -- (-.1, 1.9); 
\draw [->]  (0.15, 1.91) -- (3.7, .06);  %be
 \draw [<-]  (0.3, 2) -- (3.8, .15);
 \draw [->] (2.55, 3) -- (3.9, .15);  %db
 \draw [<-] (2.75, 2.9) -- (4.1, .15); 
 \draw [<-]  (4, 0) -- (4.9, 1.8); %bc
 \draw [->]  (4.2, 0) -- (5.05, 1.75);
\draw [->]  (4.9, 2.1) -- (2.75, 3); %cd
\draw [<-]  (4.9, 2.23) -- (3, 3.1); 
\draw [<-] (4.75, 2.07) -- (0.5, 2.07); %ce
\draw [->] (4.75, 1.95) -- (0.7, 1.95);
\draw [<-] (0.2, 2.1) -- (2.3, 3.05); %ed
\draw [->] (0.4, 2.05) -- (2.23, 2.85);

\node[above] at (2.4, -.16) {\tiny 15};  %ab  %
\node[below] at (2.4, 0.1) {\tiny 19};
\node at (4.4, 1) {\tiny 19}; %bc %
\node at (4.75, 0.75) {\tiny 23};
\node[left] at (0.5, 1) {\tiny 12}; %ae %
\node  at (0.65, .8) {\tiny 22}; 
\node at (3.75, 2.5) {\tiny 9};  %cd %
\node[above] at (3.75, 2.55) {\tiny 25}; 
\node[above] at (1.25, 2.5) {\tiny 19}; %ed %
\node  at (1.25, 2.5) {\tiny 21};

%\node[below] at (3.1, 1.3) {$0$}; %ac
\node[right] at (2.64, .89) {\tiny 26}; %AC %
\node[right] at (2.35, 1.15) {\tiny 14}; 
%\node[left] at (1.8, 1.75) {$2$}; %ad
\node[left] at (2.16, 1.75) {\tiny 21}; %ad  %
\node at (2.12, 1.55) {\tiny 23}; 
%\node[below] at (1.9, 1.3) {$1$}; %eb  %
\node[right] at (1.42, .89) {\tiny 24};
\node[right] at (1.69, 1.1) {\tiny 6};
%\node[right] at (3.2, 1.75) {$0$}; %bd
\node[above] at (3.15, 1.35) {\tiny 32}; %bd  %
\node[above] at (3.54, 1.08) {\tiny 14}; 
%\node[above] at (2.5, 2) {$0$}; %ec
\node[left] at (2.9, 1.95) {\tiny 24}; %ec
\node at (2.81, 2.11) {\tiny 0}; 

\node[below] at (2.5, -1) {{\bf a.} A general $\mathcal G^5$}; 
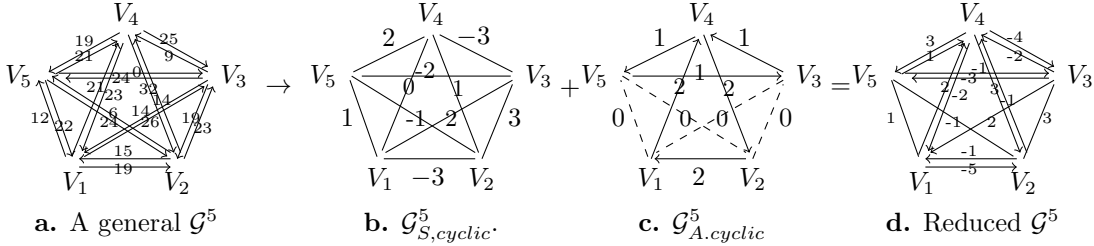
\node[below] at (12.5, -2) {{\bf Figure 14.} Decomposing an asymmetric  $\mathcal G^5$};
\node at (7, 1.7) {$\to$};

\node at (15.6, 1.7) {+};

%%%%%%%%%%%%%%%

\node[below] at (10, 0) {$V_1$};
%\node at (24, 0) {$\bullet$};
\node[below] at (13, 0) {$V_2$};

\node[right] at (14, 2) {$V_3$};
 
\node[above] at (11.5, 3) {$V_4$};
 
\node[left] at (9, 2) {$V_5$};
\draw  (10.1,  0) -- (12.8, 0); %ab
\draw   (13.9, 1.9) -- (10.15, .15);  %ca
\draw  (10.05, 0.2) -- (11.4, 2.9); %ad
\draw  (9.9, .05) -- (9.1, 1.9); %ae
\draw    (9.15, 2) -- (12.8, .15);  %be
 \draw   (11.55, 3) -- (12.9, .15);  %db
 \draw    (13, 0) -- (13.9, 1.8); %bc
\draw  (13.9, 2.1) -- (11.75, 3); %cd
\draw   (13.85, 2) -- (9.45, 2); %ce
\draw  (9.2, 2.1) -- (11.3, 2.95); %ed

\node[below] at (11.4, 0) {$-3$};  %ab
\node[right] at (13.5, 1) {$3$}; %bc
\node[left] at (9.5, 1) {$1$}; %ae
\node[above] at (12.75, 2.5) {$-3$};  %cd
\node[above] at (10.25, 2.5) {$2$}; %ed

%\node[below] at (3.1, 1.3) {$0$}; %ac
\node[right] at (11.64, .99) {\small 2}; %AC
%\node[left] at (1.8, 1.75) {$2$}; %ad
\node[left] at (11.3, 1.75) {\small 0}; %ad
%\node[below] at (1.9, 1.3) {$1$}; %eb
\node[right] at (10.48, .99) {\small -1};
%\node[right] at (3.2, 1.75) {$0$}; %bd
\node[above] at (12.3, 1.3) {\small 1}; %bd
%\node[above] at (2.5, 2) {$0$}; %ec
\node[left] at (11.9, 2.15) {\small -2}; %ec

\node[below] at (11.5, -1) {{\bf b.} $\mathcal G^5_{S, cyclic}$.}; 

%%%%%%%%%%%%%%%%%%%%%
%ab=1, ac=-4, ad=-3, ae=-9, bc=-5, bd=-4, be=-10; cd=1, ce=-5, de=-6
\node[below] at (18, 0) {$V_1$};
%\node at (24, 0) {$\bullet$};
\node[below] at (21, 0) {$V_2$};

%\node at (25, 2) {$\bullet$};
\node[right] at (22, 2) {$V_3$};
%\node at (22.5, 3) {$\bullet$};
\node[above] at (19.5, 3) {$V_4$};
%\node at (20, 2) {$\bullet$};
\node[left] at (17, 2) {$V_5$};
\draw[<-] (18.1,  0) -- (20.8, 0); %ab
\draw [dashed] (21.9, 1.9) -- (18.15, .15);  %ca
\draw [->]  (18, 0) -- (19.4, 2.9); %ad
\draw [dashed] (17.9, .05) -- (17.1, 1.9); %ae
\draw [dashed]  (17.15, 2) -- (20.8, .15);  %be
 \draw [->] (19.55, 3) -- (20.9, .15);  %db
 \draw [dashed]  (21, 0) -- (21.9, 1.8); %bc
\draw [->]  (21.9, 2.1) -- (19.75, 3); %cd
\draw [<-] (21.85, 2) -- (17.45, 2); %ce
\draw [<-] (17.2, 2.1) -- (19.3, 2.95); %ed

\node[below] at (19.4, 0) {$2$};  %ab
\node[right] at (21.5, 1) {$0$}; %bc
\node[left] at (17.5, 1) {$0$}; %ae
\node[above] at (20.75, 2.5) {$1$};  %cd
\node[above] at (18.25, 2.5) {$1$}; %ed

%\node[below] at (3.1, 1.3) {$0$}; %ac
\node[right] at (19.64, .99) {\small 0}; %AC
%\node[left] at (1.8, 1.75) {$2$}; %ad
\node[left] at (19.3, 1.75) {\small 2}; %ad
%\node[below] at (1.9, 1.3) {$1$}; %eb
\node[right] at (18.56, .99) {\small 0};
%\node[right] at (3.2, 1.75) {$0$}; %bd
\node[above] at (20.3, 1.3) {\small 2}; %bd
%\node[above] at (2.5, 2) {$0$}; %ec
\node[left] at (19.9, 2.1) {\small 1}; %ec

\node[below] at (19.5, -1) {{\bf c.}   $\mathcal G^5_{A. cyclic}$}; 

%%%%%%%%%%%%%%%%%%%%%%

\node[below] at (26, -.10) {$V_1$};
%%\node at (24, 0) {$\bullet$};
\node[below] at (29, -.10) {$V_2$};
;

%\node at (25, 2) {$\bullet$};
\node[right] at (30, 2) {$V_3$};
%\node at (22.5, 3) {$\bullet$};
\node[above] at (27.5, 3) {$V_4$};
%\node at (20, 2) {$\bullet$};
\node[left] at (25, 2) {$V_5$};
\draw[<-] (26.1,  0) -- (28.8, 0); %ab
\draw[->] (26.1, -.2) -- (28.8, -.2);
\draw [->] (29.9, 1.9) -- (26.25, .15);  %ca
%\draw [<-] (29.8, 1.75) -- (26.25, .05);
\draw [->]  (25.9, .1) -- (27.25, 2.8); %ad
\draw [<-]  (26.19, .19) -- (27.4, 2.8); 
\draw [-] (25.9, .05) -- (25.1, 1.9); %ae
%\draw [->] (25.7, .05) -- (24.9, 1.9); 
\draw [-]  (25.15, 1.91) -- (28.7, .06);  %be
% \draw [<-]  (25.3, 2) -- (28.8, .15);
 \draw [->] (27.55, 3) -- (28.9, .15);  %db
 \draw [<-] (27.75, 2.9) -- (29.1, .15); 
 \draw [-]  (29, 0) -- (29.9, 1.8); %bc
% \draw [->]  (29.2, 0) -- (30.05, 1.75);
\draw [->]  (29.9, 2.1) -- (27.75, 3); %cd
\draw [<-]  (29.9, 2.23) -- (28, 3.1); 
\draw [<-] (29.75, 2.07) -- (25.5, 2.07); %ce
\draw [->] (29.75, 1.95) -- (25.7, 1.95);
\draw [<-] (25.2, 2.1) -- (27.3, 3.05); %ed
\draw [->] (25.4, 2.05) -- (27.23, 2.85);

\node[above] at (27.4, -.16) {\tiny -1};  %ab
\node[below] at (27.4, 0.1) {\tiny -5};
\node at (29.7, 1) {\tiny 3}; %bc
%\node at (29.75, 0.75) {\tiny 15};
\node[left] at (25.5, 1) {\tiny 1}; %ae
%\node  at (25.65, .8) {\tiny 21}; 
\node at (28.75, 2.5) {\tiny -2};  %cd
\node[above] at (28.75, 2.59) {\tiny -4}; 
\node[above] at (26.25, 2.5) {\tiny 3}; %ed
\node  at (26.25, 2.5) {\tiny 1};

%\node[below] at (3.1, 1.3) {$0$}; %ac
\node[right] at (27.64, .89) {\tiny 2}; %AC
 
\node[left] at (27.1, 1.75) {\tiny 2}; %ad
\node at (27.12, 1.55) {\tiny -2}; 
%\node[below] at (1.9, 1.3) {$1$}; %eb
\node[right] at (26.42, .89) {\tiny -1};

\node[above] at (28.15, 1.35) {\tiny 3}; %bd
\node[above] at (28.54, 1.08) {\tiny -1}; 
%\node[above] at (2.5, 2) {$0$}; %ec
\node[left] at (27.9, 1.95) {\tiny -3}; %ec
\node at (27.7, 2.19) {\tiny -1}; 

\node[below] at (27.5, -1) {{\bf d.} Reduced $\mathcal G^5$}; 

\node at (23.6, 1.7) {=};

\end{tikzpicture}

To illustrate this program with Fig.~14,  
  by representing  the Fig.~14a  costs    as   \{average cost, excess cost\} pairs,  the original graph becomes  
    $\mathcal G^5 = \mathcal G^5_S +\mathcal G^5_A \in \mathbb G^5_S \times \mathbb G^5_A$, where  $\widehat{V_jV_k}$'s   length in  $\mathcal G^5_S$ is the average cost of its   arcs, and  $\mathcal G^5_A$ 
  represents
how   costs differ from the average.   Thus,  with 
   $V_1 \stackrel{26}{\longrightarrow} V_3$ and $V_3 \stackrel{14}{\longrightarrow} V_1$ from Fig.~14a,  $\widehat{V_1V_3}$'s   length in $\mathcal G^5_S$ is $20$ and  
  $\mathcal G^5_A$ has  $V_1 \stackrel{6}{\longrightarrow} V_3$.  

The analysis of $\mathcal G^5_S + \mathcal G^5_A$ follows as above:  find each graph's cpi and cyclic components.    Removing      $\mathcal G^5_{S, cpi}$ and $\mathcal G^5_{A, cpi}$  leaves   Figs.~14b, c.  
  A $\mathcal G^5$  Hamiltonian path length  
  (with   $\mathcal G^5_{S, cpi}$ weights $\omega_1=10, \,  \omega_2=10,  \,  \omega_3 = 8, \, \omega_4= 12, \, \omega_5= 6$)  equals $T(\mathcal G^5_{S, cpi}) = 92$    
  plus the sum of its   $\mathcal G^5_{S, cyclic}$ and $\mathcal G^5_{A, cyclic}$ lengths.   Expressing   $\mathcal G^5_{S, cyclic} + \mathcal G^5_{A, cyclic}$ in a  standard  Fig.~14d form,   its shortest Hamiltonian path of  
 $ V_1  \stackrel{-5}{\longrightarrow} V_2  \stackrel{-1}{\longrightarrow} V_4  \stackrel{-4}{\longrightarrow} V_3  \stackrel{-3}{\longrightarrow} V_5  \stackrel{1}{\longrightarrow} V_1$  
  follows.   
 Its $\mathcal G^5$  
 path  length  is $T(\mathcal G^5_{S, cpi}) -12 = 80$.

\section{Summary}  
Components of a graph's entries that hamper finding   closed path properties are identified.  Eliminating them determines the graph's essence---a  reduced graph with   smaller   degrees of freedom  where all  closed graph properties are expressed in terms of  inherent symmetry structures.

\section{Proofs} \label{sect: proofs}  Results not proved above or in \cite{tandd}, are proved here.

{\em Proof of Theorem 1:} For a triplet $\{V_i, V_j, V_k\}$ in a cpi graph,  the closed path  $V_i   \stackrel{x}{\longrightarrow} V_j \stackrel{y}{\longrightarrow}  V_k \stackrel{z}{\longrightarrow}  V_i$ has path length zero, so  $x+y = -z.$  Thus, $V_i   \stackrel{x}{\longrightarrow} V_j \stackrel{y}{\longrightarrow}  V_k = V_i \stackrel{-z= x+y}{\longrightarrow}  V_k$ satisfies Eq.~\ref{eq: st}.  As all     triplets satisfy Eq.~\ref{eq: st},  
   a cpi graph is strongly transitive.

A triplet $\{V_i, V_j, V_s\}$ in a strongly transitive graph satisfies  $V_i \stackrel{x}{\longrightarrow} V_j \stackrel{y}{\longrightarrow} V_s = V_i \stackrel{z}{\longrightarrow} V_s$ where $z=x+y$.  Applying  a fourth alternative $V_t$ to this relationship yields $$(V_i \stackrel{x}{\longrightarrow} V_j \stackrel{y}{\longrightarrow} V_s) \stackrel{u}{\longrightarrow} V_t =  (V_i \stackrel{z}{\longrightarrow} V_s) \stackrel{u}{\longrightarrow} V_t   = V_i \stackrel{w}{\longrightarrow} V_t,$$ where $w=z+u = x+y+u.$  With the obvious induction argument, it follows that any path from $V_i$ to $V_k$ has the same  length as the direct path from $V_i$ to $V_k$.  (This proves Cor.~1.)  A closed path has $V_k=V_i$, so its  length is that of $V_i$ to $V_i$, or zero.  Hence, a strongly transitive graph is cpi.

To prove that the set of strongly transitive graphs forms a linear subspace, notice that a multiple $\mu$ of a  
 strongly transitive graph in $\mathbb G^n_A$ changes all  path lengths  by this multiple; thus the new graph's   arcs remain  strongly transitive.  Therefore the multiple  defines another $\mathbb G^n_A$ strongly transitive graph. (If $\mu<0$, then positive cost arcs in the original graph become negative cost arcs in the new graph.)  Similarly, for two strongly transitive  $\mathbb G^n_A$ graphs and any $\{V_i, V_j, V_k\}$ triplet, the first graph satisfies  $V_i \stackrel{x}{\longrightarrow} V_j  \stackrel{y}{\longrightarrow} V_k = V_i \stackrel{x+y}{\longrightarrow} V_k$ while the second satisfies $ V_i \stackrel{\tilde{x}}{\longrightarrow} V_j  \stackrel{\tilde{y}}{\longrightarrow} V_k = V_i \stackrel{\tilde{x}+\tilde{y}}{\longrightarrow} V_k.$  Combining these graphs leads to $ V_i \stackrel{x+\tilde{x}}{\longrightarrow} V_j  \stackrel{y + \tilde{y}}{\longrightarrow} V_k = V_i \stackrel{(x+y) + (\tilde{x} + \tilde{y})}{\longrightarrow} V_k$, which satisfies Eq.~\ref{eq: st}.  
  Thus the set of strongly transitive graphs in $\mathbb G^n_A$, $\mathbb {ST}_A^n$,  is  a linear subspace.  $\square$\smallskip

{\em Proof of Cor.~2:}  The Eq.~4 basis of $ \{V_1 \stackrel{1}{\longrightarrow} V_j \stackrel{1}{\longrightarrow} V_k \stackrel{1}{\longrightarrow} V_1\}_{1<j<k\le n}$ satisfies Cor.~2  because in this set,  only the three cycle $ V_1 \stackrel{1}{\longrightarrow} V_s \stackrel{1}{\longrightarrow} V_k \stackrel{1}{\longrightarrow} V_1$ has a $\widehat{V_sV_k}$ arc. 
For the  independence of the arcs, if $\mathbf c_{j,k}$ represents the only  $\mathcal{CB}^n_A$ cycle with a $\widehat{V_jV_k}$ arc,   it must be shown that $\sum x_{j,k}\mathbf c_{j,k} = \mathbf 0$ iff all $x_{j,k}=0.$  But as $\mathbf c_{j,k}$ is the only vector with a non-zero ${j, k}$ component, $x_{j,k}=0.$   

That theses cycles are in $\mathbb{ST}_A^n$'s normal bundle is proved in \cite{tandd}.  
As this set consists of ${n-1}\choose2$ linearly independent elements that are orthogonal to $\mathbb{ST}^n_A$, it is a basis for the normal bundle.  $\square$\smallskip

{\em Proof of Thm.~3:}  Equation~\ref{eq: splitting} is an immediate consequence of the  representation of $\mathbb G_A^n$ into the orthogonal subspaces $\mathbb{ST}_A^n$ and $\mathbb C^n_A$. The last comment is proved above.$\square$\smallskip

{\em Proof of Cor.~5:} That the Def.~3 relationship is an equivalence relationship (reflexive, symmetricaA, transitive) follows immediately from the equality of the cyclic components.  The difference between any two graphs in $\mathbb G^n_A$ is the difference between their cpi and cyclic components.  As their cyclic components agree, the difference is the difference between cpi components.  Because $\mathbb{ST}_A^n$ is a linear subspace, this difference also is in $\mathbb{ST}_A^n$. $\square$\smallskip

{\em Proof of Thm.~6:}  The fact that $\mathcal G^n_{A, cpi}$ normally has a sink and source for positive and negative directions follows from the fact that all triplets are transitive, so there is a maximum and a minimum term.   That its arc lengths are non-zero means that this top and bottom alternative are unique.  For positive directions, the top alternative  is a source, the bottom alternative is a sink.  

The assertion that  $\mathcal G^n_{A, cyclic}$ cannot have a source or a sink follows from the fact that $S_A(V_j)=0$ for each vertex.  (This statement follows from the fact that each three-cycle attache to a vertex has one leg pointing in and one leg, of same magnitude, pointing out.)  Thus, each  $\mathcal G^n_{A, cyclic}$  vertex with non-zero arcs has at least one positive direction pointing in and at least one pointing out.

What remains is to show that the longest Hamiltonian path in $\mathcal G^n_{A, cyclic}$, $n=4, 5$, has all positive directions.  It already has been shown that $\mathcal G^n_{A, cyclic}$ does not have a sink or source.  The next possible problem is three-cycle with, say, positive cost arcs, e.g., $V_2\to V_3 \to V_4 \to V_2$.  To avoid having all positive arcs in the longest paths, the cycle must be attracting.   But for $n=4$, that would require all positive cost directions to point away from $V_1$, making $V_1$ a source, which it cannot be.  For $n=5$, all  positive cost arrows from $V_1$ and $V_5$ point to the cycle.  The positive cost arrow between $V_1$ and $V_5$ points away from   one of these vertices, making it a source, which is a contradiction.  An attracting four cycle would make the remaining vertex a source. Everything extends in the same manner for sinks and for negative cost directions.    $\square$\smallskip

{\em Proof of Thm.~8:}  While the  linear algebra  proof used for $n=4$  extends,  
  an iterative argument  provides  insight.  One set of closed paths involves the vertices $\{V_1, V_2, V_3, V_4\}$; these cpi 
   graphs are based on weights $\omega_j$ attached to $V_j$, $j=1, \dots, 4$. 
Increasing the graph size to involve vertices $\{V_1, V_2, V_3, V_4, V_5\}$  requires analyzing all closed graphs in  $\{V_1, V_2, V_3,  V_5\}$.  Again, the solution  has weights $\omega'_k, \, k=1, 2, 3, 5$ assigned to the appropriate vertices.  To be consistent with arc lengths in the first 4-tuple, it must be that $\omega'_k=\omega_k$ for $k=1, 2, 3$.  Continuing in this iterative manner extends  
  the proof  to   all vertices. That the length of
 a closed path passing through the vertices $\{V_j\}_{j\in \mathcal D}$ once is $2\sum_{j\in {\mathcal D}} \omega_j$ is an immediate computation.  As a path enters and leaves vertex $V_j$, the length is increased by $2\omega_j.$   $\square$\smallskip

{\em Proof of Thm.~9:}  Set  $\{\mathbf B^n_j\}_{j=1}^{n-1}$ is independent    because only $\mathbf B^n_j$ has a non-zero $d_{j, n}$ coordinate.  If $\{\mathbf B^n_j\}_{j=1}^{n}$ is not independent,   there is a summation $\sum_{j=1}^{n-1} x_j\mathbf B^n_j = \mathbf B^n_n.$  In the sum, each $x_j=1$ to capture $\mathbf B^n_n$'s $d_{j, n}=1$  component.  But then  $d_{1, 2}=2$   (from $\mathbf B^n_1 + \mathbf B^n_2$), rather than the required zero of $\mathbf B^n_n$,  so the linear subspace spanned by $\{\mathbf B^n_j\}_{j=1}^{n}$ is $n$-dimensional.  
This space captures the structure of   $\mathcal G^n_{S, cpi}$ graphs  because the $d_{i, j}$ component of  $\sum_{s=1}^n \omega_s \mathbf B^n_s$ is the required $\omega_i+\omega_j$. $\square$   

{\em Proof of Thm.~14:}   
A path's length   in $\mathcal G^n_S$ is the sum of its lengths  in $\mathcal G^n_{S, cpi}$ and in $\mathcal G^n_{S, cyclic}$.  $\square$

\end{document}